\magnification1200
\newread\AUX\immediate\openin\AUX=\jobname.aux
\def\ref#1{\expandafter\edef\csname#1\endcsname}
\ifeof\AUX\immediate\write16{\jobname.aux gibt es nicht!}\else
\input \jobname.aux
\fi\immediate\closein\AUX
\def\today{\number\day.~\ifcase\month\or
  Januar\or Februar\or M{\"a}rz\or April\or Mai\or Juni\or
  Juli\or August\or September\or Oktober\or November\or Dezember\fi
  \space\number\year}
\font\sevenex=cmex7
\scriptfont3=\sevenex
\font\fiveex=cmex10 scaled 500
\scriptscriptfont3=\fiveex

\def\G{{\bf G}}

\def\phi{\varphi}
\def\epsilon{\varepsilon}
\def\theta{\vartheta}
\def\uauf{\lower1.7pt\hbox to 3pt{%
\vbox{\offinterlineskip
\hbox{\vbox to 8.5pt{\leaders\vrule width0.2pt\vfill}%
\kern-.3pt\hbox{\lams\char"76}\kern-0.3pt%
$\raise1pt\hbox{\lams\char"76}$}}\hfil}}
\def\cite#1{\expandafter\ifx\csname#1\endcsname\relax
{\bf?}\immediate\write16{#1 ist nicht definiert!}\else\csname#1\endcsname\fi}
\def\expandwrite#1#2{\edef\next{\write#1{#2}}\next}
\def\neverexpand{\noexpand\noexpand\noexpand}
\def\strip#1\ {}
\def\ncite#1{\expandafter\ifx\csname#1\endcsname\relax
{\bf?}\immediate\write16{#1 ist nicht definiert!}\else
\expandafter\expandafter\expandafter\strip\csname#1\endcsname\fi}
\newwrite\AUX
\immediate\openout\AUX=\jobname.aux
\newcount\Abschnitt\Abschnitt0
\def\beginsection#1. #2 \par{\advance\Abschnitt1%
\vskip0pt plus.10\vsize\penalty-250
\vskip0pt plus-.10\vsize\bigskip\vskip\parskip
\edef\TEST{\number\Abschnitt}
\expandafter\ifx\csname#1\endcsname\TEST\relax\else
\immediate\write16{#1 hat sich geaendert!}\fi
\expandwrite\AUX{\neverexpand\ref{#1}{\TEST}}
\leftline{\bf\number\Abschnitt. \ignorespaces#2}%
\nobreak\smallskip\noindent\SATZ1}
\def\Proof:{\par\noindent{\it Proof:}}
\def\Remark:{\ifdim\lastskip<\medskipamount\removelastskip\medskip\fi
\noindent{\bf Remark:}}
\def\Remarks:{\ifdim\lastskip<\medskipamount\removelastskip\medskip\fi
\noindent{\bf Remarks:}}
\def\Definition:{\ifdim\lastskip<\medskipamount\removelastskip\medskip\fi
\noindent{\bf Definition:}}
\def\Example:{\ifdim\lastskip<\medskipamount\removelastskip\medskip\fi
\noindent{\bf Example:}}
\def\Examples:{\ifdim\lastskip<\medskipamount\removelastskip\medskip\fi
\noindent{\bf Examples:}}
\newcount\SATZ\SATZ1
\def\proclaim #1. #2\par{\ifdim\lastskip<\medskipamount\removelastskip
\medskip\fi
\noindent{\bf#1.\ }{\it#2}\Par
\ifdim\lastskip<\medskipamount\removelastskip\goodbreak\medskip\fi}
\def\Aussage#1{%
\expandafter\def\csname#1\endcsname##1.{\ifx?##1?\relax\else
\edef\TEST{#1\penalty10000\ \number\Abschnitt.\number\SATZ}
\expandafter\ifx\csname##1\endcsname\TEST\relax\else
\immediate\write16{##1 hat sich geaendert!}\fi
\expandwrite\AUX{\neverexpand\ref{##1}{\TEST}}\fi
\proclaim {\number\Abschnitt.\number\SATZ. #1\global\advance\SATZ1}.}}
\Aussage{Theorem}
\Aussage{Proposition}
\Aussage{Corollary}
\Aussage{Lemma}
\font\la=lasy10
\def\strich{\hbox{$\vcenter{\hbox
to 1pt{\leaders\hrule height -0,2pt depth 0,6pt\hfil}}$}}
\def\dashedrightarrow{\hbox{%
\hbox to 0,5cm{\leaders\hbox to 2pt{\hfil\strich\hfil}\hfil}%
\kern-2pt\hbox{\la\char\string"29}}}

\def\Bindestrich{\penalty10000-\hskip0pt}
\let\str=\_
\let\_=\Bindestrich
\def\.{{\sfcode`.=1000.}}

\def\Par{\par}
\def\:={\mathrel{\raise0,9pt\hbox{.}\kern-2,77779pt
\raise3pt\hbox{.}\kern-2,5pt=}}
\def\=:{\mathrel{=\kern-2,5pt\raise0,9pt\hbox{.}\kern-2,77779pt
\raise3pt\hbox{.}}} 
\def\into{\hookrightarrow}
\def\pfeil{\rightarrow}

\def\Ugleich{\hbox{$\cup$\kern.5pt\vrule depth -0.5pt}}
\def\|#1|{\mathop{\rm#1}\nolimits}
\def\<{\langle}
\def\>{\rangle}
\let\Times=\times
\def\times{\mathop{\Times}}
\let\Otimes=\otimes
\def\otimes{\mathop{\Otimes}}
\catcode`\@=11
\def\hex#1{\ifcase#1 0\or1\or2\or3\or4\or5\or6\or7\or8\or9\or A\or B\or
C\or D\or E\or F\else\message{Warnung: Setze hex#1=0}0\fi}
\def\fontdef#1:#2,#3,#4.{%
\alloc@8\fam\chardef\sixt@@n\FAM
\ifx!#2!\else\expandafter\font\csname text#1\endcsname=#2
\textfont\the\FAM=\csname text#1\endcsname\fi
\ifx!#3!\else\expandafter\font\csname script#1\endcsname=#3
\scriptfont\the\FAM=\csname script#1\endcsname\fi
\ifx!#4!\else\expandafter\font\csname scriptscript#1\endcsname=#4
\scriptscriptfont\the\FAM=\csname scriptscript#1\endcsname\fi
\expandafter\edef\csname #1\endcsname{\fam\the\FAM\csname text#1\endcsname}
\expandafter\edef\csname hex#1fam\endcsname{\hex\FAM}}
\catcode`\@=12 

\fontdef Ss:cmss10,,.
\fontdef Fr:eufm10,eufm7,eufm5.


\def\ft{{\Fr t}}

\fontdef bbb:msbm10,msbm7,msbm5.
\fontdef mbf:cmmib10,cmmib7,.
\fontdef msa:msam10,msam7,msam5.
\def\CC{{\bbb C}}

\def\NN{{\bbb N}}\def\OO{{\bbb O}}
\def\QQ{{\bbb Q}}

\def\ZZ{{\bbb Z}}
\def\cA{{\cal A}}\def\cD{{\cal D}}
\def\cH{{\cal H}}

\def\cO{{\cal O}}\def\cP{{\cal P}}

\mathchardef\leer=\string"0\hexbbbfam3F
\mathchardef\subsetneq=\string"3\hexbbbfam24
\mathchardef\semidir=\string"2\hexbbbfam6E
\mathchardef\dirsemi=\string"2\hexbbbfam6F
\mathchardef\haken=\string"2\hexmsafam78
\mathchardef\auf=\string"3\hexmsafam10
\let\OL=\overline
\def\overline#1{{\hskip1pt\OL{\hskip-1pt#1\hskip-1pt}\hskip1pt}}


\def\Gq{{\overline{G}}}

\def\Uq{{\overline{U}}}

\def\Wq{{\overline{W}}}

\def\Zq{{\overline{Z}}}
%
\abovedisplayskip 9.0pt plus 3.0pt minus 3.0pt
\belowdisplayskip 9.0pt plus 3.0pt minus 3.0pt
\newdimen\Grenze\Grenze2\parindent\advance\Grenze1em
\newdimen\Breite
\newbox\DpBox
\def\NewDisplay#1$${\Breite\hsize\advance\Breite-\hangindent
\setbox\DpBox=\hbox{\hskip2\parindent$\displaystyle{#1}$}%
\ifnum\predisplaysize<\Grenze\abovedisplayskip\abovedisplayshortskip
\belowdisplayskip\belowdisplayshortskip\fi
\global\futurelet\nexttok\WEITER}
\def\WEITER{\ifx\nexttok\qed\expandafter\leftQEDdisplay
\else\leftdisplay\fi}
\def\leftdisplay{\hskip-\hangindent\leftline{\box\DpBox}$$}
\def\leftQEDdisplay{\hskip-\hangindent
\line{\copy\DpBox\hfill\lower\dp\DpBox\copy\QEDbox}%
\belowdisplayskip0pt$$\bigskip\let\nexttok=}
\everydisplay{\NewDisplay}
\newbox\QEDbox
\newbox\nichts\setbox\nichts=\vbox{}\wd\nichts=2mm\ht\nichts=2mm
\setbox\QEDbox=\hbox{\vrule\vbox{\hrule\copy\nichts\hrule}\vrule}
\def\qed{\leavevmode\unskip\hfil\null\nobreak\hfill\copy\QEDbox\medbreak}
\newdimen\HIindent
\newbox\HIbox
\def\setHI#1{\setbox\HIbox=\hbox{#1}\HIindent=\wd\HIbox}
\def\HI#1{\par\hangindent\HIindent\hangafter=0\noindent\leavevmode
\llap{\hbox to\HIindent{#1\hfil}}\ignorespaces}

\newdimen\maxSpalbr
\newdimen\altSpalbr
\newcount\Zaehler

\def\beginrefs{%
\expandafter\ifx\csname Spaltenbreite\endcsname\relax
\def\Spaltenbreite{1cm}\immediate\write16{Spaltenbreite undefiniert!}\fi
\expandafter\altSpalbr\Spaltenbreite
\maxSpalbr0pt
\gdef\alt{}
\def\\##1\relax{%
\gdef\neu{##1}\ifx\alt\neu\global\advance\Zaehler1\else
\xdef\alt{\neu}\global\Zaehler=1\fi\xdef\SigText{##1\the\Zaehler}}
\def\L|Abk:##1|Sig:##2|Au:##3|Tit:##4|Zs:##5|Bd:##6|S:##7|J:##8||{%
\def\SigText{##2}\global\setbox0=\hbox{##2\relax}
\edef\TEST{[\SigText]}
\expandafter\ifx\csname##1\endcsname\TEST\relax\else
\immediate\write16{##1 hat sich geaendert!}\fi
\expandwrite\AUX{\neverexpand\ref{##1}{\TEST}}
\setHI{[\SigText]\ }
\ifnum\HIindent>\maxSpalbr\maxSpalbr\HIindent\fi
\ifnum\HIindent<\altSpalbr\HIindent\altSpalbr\fi
\HI{[\SigText]}
\ifx-##3\relax\else{##3}: \fi
\ifx-##4\relax\else{##4}{\sfcode`.=3000.} \fi
\ifx-##5\relax\else{\it ##5\/} \fi
\ifx-##6\relax\else{\bf ##6} \fi
\ifx-##8\relax\else({##8})\fi
\ifx-##7\relax\else, {##7}\fi\Par}
\def\B|Abk:##1|Sig:##2|Au:##3|Tit:##4|Reihe:##5|Verlag:##6|Ort:##7|J:##8||{%
\def\SigText{##2}\global\setbox0=\hbox{##2\relax}
\edef\TEST{[\SigText]}
\expandafter\ifx\csname##1\endcsname\TEST\relax\else
\immediate\write16{##1 hat sich geaendert!}\fi
\expandwrite\AUX{\neverexpand\ref{##1}{\TEST}}
\setHI{[\SigText]\ }
\ifnum\HIindent>\maxSpalbr\maxSpalbr\HIindent\fi
\ifnum\HIindent<\altSpalbr\HIindent\altSpalbr\fi
\HI{[\SigText]}
\ifx-##3\relax\else{##3}: \fi
\ifx-##4\relax\else{##4}{\sfcode`.=3000.} \fi
\ifx-##5\relax\else{(##5)} \fi
\ifx-##7\relax\else{##7:} \fi
\ifx-##6\relax\else{##6}\fi
\ifx-##8\relax\else{ ##8}\fi\Par}
\def\Pr|Abk:##1|Sig:##2|Au:##3|Artikel:##4|Titel:##5|Hgr:##6|Reihe:{%
\def\SigText{##2}\global\setbox0=\hbox{##2\relax}
\edef\TEST{[\SigText]}
\expandafter\ifx\csname##1\endcsname\TEST\relax\else
\immediate\write16{##1 hat sich geaendert!}\fi
\expandwrite\AUX{\neverexpand\ref{##1}{\TEST}}
\setHI{[\SigText]\ }
\ifnum\HIindent>\maxSpalbr\maxSpalbr\HIindent\fi
\ifnum\HIindent<\altSpalbr\HIindent\altSpalbr\fi
\HI{[\SigText]}
\ifx-##3\relax\else{##3}: \fi
\ifx-##4\relax\else{##4}{\sfcode`.=3000.} \fi
\ifx-##5\relax\else{In: \it ##5}. \fi
\ifx-##6\relax\else{(##6)} \fi\PrII}
\def\PrII##1|Bd:##2|Verlag:##3|Ort:##4|S:##5|J:##6||{%
\ifx-##1\relax\else{##1} \fi
\ifx-##2\relax\else{\bf ##2}, \fi
\ifx-##4\relax\else{##4:} \fi
\ifx-##3\relax\else{##3} \fi
\ifx-##6\relax\else{##6}\fi
\ifx-##5\relax\else{, ##5}\fi\Par}
\bgroup
\baselineskip12pt
\parskip2.5pt plus 1pt
\hyphenation{Hei-del-berg}
\sfcode`.=1000
\beginsection References. References

}
\def\endrefs{%
\expandwrite\AUX{\neverexpand\ref{Spaltenbreite}{\the\maxSpalbr}}
\ifnum\maxSpalbr=\altSpalbr\relax\else
\immediate\write16{Spaltenbreite hat sich geaendert!}\fi
\egroup}


\fontdef Sans:cmss10,,.

\def\da{\downarrow}
\def\1{{\textstyle{1\over2}}}

\newcount\GNo\GNo=0
\def\eqno#1{
\global\advance\GNo1
\edef\FTEST{$(\number\GNo)$}
\ifx?#1?\relax\else
\expandafter\ifx\csname#1\endcsname\FTEST\relax\else
\immediate\write16{#1 hat sich geaendert!}\fi
\expandwrite\AUX{\neverexpand\ref{#1}{\FTEST}}\fi
\llap{\hbox to 40pt{\FTEST\hfill}}}

\newcount\Condition
\Condition=0
\def\cond#1.{
\global\advance\Condition1
\edef\FTEST{{\bf C\the\Condition}}
\ifx?#1?\relax\else
\expandafter\xdef\csname Condition#1\endcsname{{\bf C\the\Condition}}
\fi
\llap{\hbox to 40pt{\FTEST\hfill}}}

\def\co#1.{\expandafter
\ifx\csname
Condition#1\endcsname\relax\immediate\write16{C#1 ist nicht
definiert!}{\bf C?}\else\csname Condition#1\endcsname\fi}

\font\BF=cmbx10 scaled \magstep2
\font\CSC=cmcsc10 
\baselineskip14pt

{\baselineskip1.5\baselineskip\rightskip0pt plus 5truecm
\leavevmode\vskip0truecm\noindent
\BF Construction of Commuting Difference Operators for Multiplicity
Free Spaces

}
\vskip1truecm
\leftline{{\CSC Friedrich Knop}*%
\footnote{}{\sevenrm \llap{*} Supported by a grant of the NSF.}}
\leftline{Department of Mathematics, Rutgers University, New Brunswick NJ
08903, USA}
\leftline{knop@math.rutgers.edu}
\vskip1truecm

\beginsection Introduction. Introduction

The analysis of invariant differential operators on certain
multiplicity free spaces led recently to the introduction of a family
of symmetric polynomials that is more general than Jack polynomials
(see \cite{SymCap}, but also \cite{OO1}, \cite{OO2}). They are called
{\it interpolation Jack polynomials}, {\it shifted Jack polynomials},
or {\it Capelli polynomials}. Apart from being inhomogeneous, they are
distinguished from classical Jack polynomials by their very simple
definition in terms of vanishing conditions.

One of the most important and non\_obvious properties of Capelli
polynomials is that they are eigenfunctions of certain explicitly
given difference (as opposed to differential) operators (see
\cite{SymCap}). This readily implies that their top homogeneous term
is in fact a (classical) Jack polynomial. Other consequences include a
binomial theorem, a Pieri formula, and much more.

It is well\_known that Jack polynomials are tied to root systems of
type $\Sans A$ and that they have natural analogues for other root
systems (see e.g., \cite{He}). Therefore, it is a natural problem as to
whether this holds for the Capelli polynomials as well. Okounkov
\cite{Ok2} proposed such an analogue for root systems of type $\Sans
BC$ and proved that these share some of the nice properties of Capelli
polynomials. But, unfortunately, Okounkov's polynomials do not satisfy
difference equations. Also their representation\_theoretic significance
is not clear.

In this paper we go back to the origin and let ourselves be guided by the
theory of multiplicity free actions. It is known (see
section~\cite{multiplicity free} for details) that these actions give rise to
combinatorial structures consisting of four data
$(\Gamma,\Sigma,W,\ell)$. Here $\Gamma$ is a lattice,
$\Sigma\subset\Gamma$ a basis, $W\subseteq\|Aut|\Gamma$ a finite
reflection group, and $\ell\in\Gamma$ some element. These data alone
suffice to {\it formulate\/} the definition of a (generalized) Capelli
polynomial but, in that generality, neither existence nor uniqueness
will hold, let alone any other good property.

It is known that for structures coming from multiplicity free actions,
Capelli polynomials {\it are\/} well\_defined. Moreover, two of the
most important cases (the classical,cite{SymCap} and the
semiclassical \cite{Semisym}, see tables in section~\cite{Tables})
have previously been worked out in detail. In particular, it was shown
that the corresponding Capelli polynomials are eigenfunctions of
certain difference operators.

In this paper we handle all other multiplicity free actions in an
axiomatic fashion. We extracted from the multiplicity free case nine
properties {\bf C1}--{\bf C9} and use them as the foundation of the
theory. The main result is the construction of a commuting family of
difference operators that is diagonalized by the Capelli
polynomials. In a forthcoming paper, we study the algebra of
difference operators in more detail and derive an evaluation formula,
an explicit interpolation formula, and a binomial formula, among
others.

It should be mentioned that the actual verification of properties {\bf
C1}--{\bf C9} requires some case\_by\_case analysis which uses the
classification results of \cite{Kac}, Benson\_Ratcliff \cite{BenRat1},
and Leahy \cite{Leahy}. On the other side, this disadvantage is offset
by two things: first, other structures $(\Gamma,\Sigma,W,\ell)$ which
do not come from multiplicity free actions may and do satisfy the
axioms. Thus the theory developed in this paper has applications
beyond multiplicity free actions even though the exact scope is not
yet known.

Secondly, as kind of a byproduct, all Capelli polynomials depend on at
least one free parameter. The polynomials that are actually attached
to a multiplicity free space correspond to one particular choice of
the parameters. This extra generality does not seem possible when
working directly with the multiplicity free action.

There is another paper, \cite{BenRat2} by Benson and Ratcliff, which
studies eigenvalues of invariant differential operators on
multiplicity free spaces from a combinatorial point of view. It is
just opposite in its approach: multiplicity free spaces are treated
conceptually but there are no parameters. Moreover, only the special
values $p_\mu(\rho+\lambda)$ are investigated and not their
interpolation $p_\mu(z)$. Nevertheless, the influence of
\cite{BenRat2} on the present treatment is acknowledged. This holds in
particular for formula \cite{ellPieri2}, apparently due to Yan
\cite{Yan}, and the realization of how much can be deduced from it.

Finally, one important difference from the Jack case should be
mentioned. This paper does {\it not\/} achieve the goal to define a
shifted version of generalized Jacobi polynomials in the sense of
Heckman \cite{He} (i.e., analogues of Jack polynomials for other root
systems): in general, the top homogeneous components of Capelli
polynomials are new.  Nevertheless, these components share a lot of
properties with Jacobi polynomials such as being eigenfunctions for
certain commuting differential operators. A unifying concept would be
very desirable.

\beginsection Data. Data and axioms*\footnote{}{%
\baselineskip8.4pt\sevenrm\llap{*}
At first glance, these data and axioms
might not seem very natural. Therefore, the
reader may wish to look first at sections \cite{multiplicity free} and
\cite{Tables} for motivational background and examples.\Par}

The goal of this paper is to study special polynomials that are
constructed from the following data:

\item{$\bullet$}a lattice $\Gamma$ of finite rank;
\item{$\bullet$}a basis $\Sigma$ of $\Gamma$;
\item{$\bullet$}a finite group $W\subseteq\|Aut|(\Gamma)$;
\item{$\bullet$}an element $\ell\in\Gamma$.

\noindent Let $V:=\|Hom|(\Gamma,\CC)$ and 
let $\cP=S^\bullet(\Gamma\otimes\CC)$ denote the polynomial functions
on $V$. The dual lattice $\Gamma^\vee$ sits inside $V$. Let
$\Sigma^\vee\subseteq\Gamma^\vee$ be the basis dual to $\Sigma$ and
$\Lambda_+:=\sum\NN\Sigma^\vee$ the monoid generated by it.

The structure $(\Gamma,\Sigma,W,\ell)$ is subject to the following
conditions:
$$
\cond 11.
\hbox{\it The group $W$ is generated by reflections on $V$.}
$$
Thus, $W$ gives rise to a unique root system
$\Delta\subseteq\Gamma$ such that all roots are primitive vectors. Let
$\Delta^\vee\subseteq\Gamma^\vee$ be the set of coroots. Let
$\Delta^+:=\{\alpha\in\Delta\mid\alpha(\Sigma^\vee)\ge0\}$.
$$\cond 1. \Delta=\Delta^+\cup(-\Delta^+)$$
In other words, $\Delta^+$ is a
system of positive roots and all elements of $\Sigma^\vee$ are dominant
with respect to it. Let $\Phi:=W\Sigma$ and
$\Phi^+:=\{\omega\in\Phi\mid\omega(\Sigma^\vee)\ge0\}$.
$$\cond 2. 
  \Phi\subseteq\Phi^+\cup(-\Phi^+).$$
$$\cond 3''. \ell\in\Gamma^W.$$
$$\cond 3. \sum\Phi^+-\sum\Delta^+=\ell.$$
$$\cond 3'. \hbox{\it$\ell(\eta)>0$ for all $\eta\in\Sigma^\vee$.}$$
\noindent
Let $\Sigma^\vee_1:=\{\gamma\in\Sigma^\vee\mid\ell(\gamma)=1\}$.
$$\cond 4. \hbox{\it If
  $\eta\in\Sigma^\vee_1$ and $\omega\in\Delta\cup\Phi$ then
  $\omega(\eta)\in\{-1,0,1\}$.}$$
$$\cond 5. \hbox{\it Any linear $W$\_invariant on $V$ is uniquely
determined by its restriction to $\Sigma_1^\vee$.}
$$
Let $\pm W$ be the group generated by $W$ and $-1$. Then we define
$$
V_0:=\{\rho\in V\mid\hbox{For all $\omega_1,\omega_2\in\Sigma$
with $\omega_1\in \pm W\omega_2$ holds $\omega_1(\rho)=\omega_2(\rho)$}\}
$$
Thus, for $\rho\in V_0$ and for every $\omega\in\Phi\cup(-\Phi)$ we can define
$k_\omega:=\omega_1(\rho)$ where $\omega_1\in
\pm W\omega\cap\Sigma$. In particular we have $k_\omega=k_{-\omega}$
for all $\omega\in\Phi$.  Recall that $\rho$ is called {\it dominant}
(resp. {\it regular}) if $\alpha(\rho)\not\in\ZZ_{<0}$
(resp. $\alpha(\rho)\ne0$) for all $\alpha\in\Delta^+$. The last axiom
is:
\newdimen\nndimen
\nndimen\hsize
\advance\nndimen-4\parindent
\setbox0=\vtop{\hsize\nndimen\noindent\it
There is a regular dominant $\rho\in V_0$ with the
following property: for every $\lambda\in\Lambda_+$ there is a unique
polynomial $p\in\cP^W$ of degree $\ell(\lambda)$ such that
$p(\rho+\mu)=\delta_{\lambda\mu}$ (Kronecker delta)
for all $\mu\in\Lambda_+$ with $\ell(\mu)\le\ell(\lambda)$.}
$$\cond 0. \box0$$
\noindent We will show (\cite{Main-Prop-E}) that these polynomials then
exist in fact for every dominant $\rho\in V_0$.

\beginsection Euler. The difference Euler operator

For any $\eta\in\Lambda_1:=W\Sigma_1^\vee$ define the
following rational function on $V$:
$$
f_\eta(z):={\prod\limits_{\omega\in\Phi:\omega(\eta)>0}(\omega(z)-k_\omega)
\over
\prod\limits_{\alpha\in\Delta:\alpha(\eta)>0}\alpha(z)}.
$$
For any $\eta\in\Gamma$ we have the shift operator $T_\eta$ on $\cP$
defined by $(T_\eta f)(z)=f(z-\eta)$. Then we can define the difference
operator $L:=\sum_{\eta\in\Lambda_1}f_\eta(z)T_\eta$.

\Proposition PropL1. The operator $L$ preserves the space of
$W$\_invariant polynomials.

\Proof: Let $f\in\cP^W$. From $\rho\in V_0$ it follows that $L(f)$ is
a $W$\_invariant rational function. By \co 1., the ideal of
$W$\_skew\_invariants is generated by
$\delta=\prod_{\alpha\in\Delta^+}\alpha$. It follows from the
definition of $L$ that $\delta L(f)$ is a skew\_invariant
polynomial. Thus $L(f)$ is a $W$\_invariant polynomial.\qed

If $\rho$ is regular dominant, then $\alpha(\rho+\lambda)\ne0$ for all
$\lambda\in\Lambda_+$. Thus, $f_\eta(\rho+\lambda)$ is defined. The
main property of $f_\eta$ is the following cut\_off property:

\Lemma CutOff. Assume $\rho$ is regular dominant. Let
$\eta\in\Lambda_1$ and $\lambda\in\Lambda_+$ with
$\mu:=\lambda-\eta\not\in\Lambda_+$. Then $f_\eta(\rho+\lambda)=0$.

\Proof: Since $\mu\not\in\Lambda_+$ there is $\omega\in\Sigma$ with
$\omega(\mu)=\omega(\lambda)-\omega(\eta)<0$. We have
$\omega(\lambda)\ge0$ because $\lambda\in\Lambda_+$. Thus,
$\omega(\eta)>0$ and therefore $\omega(\eta)=1$, by \co 4.. This
implies $\omega(\lambda)=0$. Therefore, the factor
$\omega(z)-k_\omega$ of $f_\eta$ vanishes in $z=\rho+\lambda$.\qed

\noindent This has the following consequence:

\Corollary PropL2. Assume $\rho$ to be regular dominant. For every
$\lambda\in\Lambda_+$ let $M_\lambda$ be the space of functions
$f\in\cP^W$ with $f(\rho+\mu)=0$ for every $\mu\in\Lambda^+$ with
$\ell(\mu)\le\ell(\lambda)$ and $\mu\ne\lambda$. Then
$L(M_\lambda)\subseteq M_\lambda$.

\noindent Now we define the {\it difference Euler operator\/} as
$E:=\ell-L$. Clearly, it inherits the properties expressed in the last
two propositions from $L$. Additionally:

\Proposition PropE. If $f\in\cP^W$ then $\|deg|E(f)\le\|deg|f$.

\Proof: Let $\eta\in\Sigma^\vee_1$. Then $\omega\in\Phi$,
$\omega(\eta)>0$ implies $\omega\in\Phi^+$ (from \co 2.). Similarly for
$\Delta$ (from \co 1.). Since, by \co 3., we have
$\sum_{\omega\in\Phi^+}\omega(\eta)=
1+\sum_{\alpha\in\Delta^+}\alpha(\eta)$ we conclude (from \co 4.) that
there is one more $\omega\in\Phi$ with $\omega(\eta)>0$ than
$\alpha\in\Delta$ with $\alpha(\eta)>0$. Thus $f_\eta$ is a rational
function of degree one. By $W$\_equivariance, the same holds for all
$\eta\in\Lambda_1$. This shows $\|deg|L(f)\le\|deg|f+1$.

We examine the effect of $L$ on the highest degree component of
$f$. There, $T_\eta$ acts as identity. Hence, $L$ acts as
multiplication by a $W$\_invariant linear function $\ell'(z)$ which is
independent of $\rho$. It remains to be shown that $\ell'=\ell$.

Let $\eta_1\in\Sigma^\vee_1$. We enumerate the other elements of
$\Sigma^\vee$ as $\eta_2,\ldots,\eta_r$. Let $\omega_i\in\Sigma$ such
that $\omega_i(\eta_j)=\delta_{ij}$. When we write
$\ell=\sum_ia_i\omega_i$, then $a_i=\ell(\eta_i)>0$ by \co
3'.. Substituting any $\eta\in\Lambda_1$ we obtain
$$
1=\ell(\eta)=\omega_1(\eta)+\sum_{i\ge2}a_i\omega_i(\eta)\le1+
\sum_{i\ge2}a_i\omega_i(\eta).
$$
Assume $\omega_i(\eta)\le0$ for all $i\ge2$. Then
$\omega_i(\eta)=0$ for all $i\ge2$ which implies that $\eta$ is a
multiple of $\eta_1$. Because $\ell(\eta)=\ell(\eta_1)=1$ we get
$\eta=\eta_1$. Thus, for $\eta\ne\eta_1$ there is $i\ge2$ with
$\omega_i(\eta)>0$. Therefore, the factor $\omega_i(z)-k_{\omega_i}$
appears in the definition of $f_\eta$, which implies that
$f_\eta(\rho+\eta_1)=0$.

Now we apply $L$ to $f=1$ and obtain
$$\eqno{Formula3}
\sum_\eta
f_\eta(z)=L(1)=\ell'(z)+a(\rho)
$$
where $a(\rho)$ is some constant. For $\rho=0$ the left-hand side of
\cite{Formula3} becomes homogeneous, which implies $a(0)=0$. But first
we substitute $z=\rho+\eta_1$ in \cite{Formula3} where
$\eta_1\in\Sigma^\vee_1$ and get
$f_{\eta_1}(\rho+\eta_1)=\ell'(\rho+\eta_1)+a(\rho)$. Now we put
$\rho=0$ and get $1=\ell'(\eta_1)$ (by \co 4.). Thus
$\ell(\eta_1)=\ell'(\eta_1)$ for all $\eta_1\in\Sigma^\vee_1$. Since both
$\ell$ and $\ell'$ are $W$\_invariant, \co 5. implies $\ell=\ell'$.\qed

\Lemma Diagonal. The action of $E$ on $\cP^W$ is
diagonalizable. Moreover, if $\rho$ is dominant and $g\in\cP^W$ is an
eigenvector of $E$, then its eigenvalue equals
$$\eqno{EigValE}
\ell(\rho)+\|min|\{\ell(\lambda)\mid \lambda\in\Lambda_+,
g(\rho+\lambda)\ne0\}.
$$

\Proof: Assume first that $\rho$ regular and dominant. For $d\in\NN$
let $U_d$ be the space of $g\in\cP^W$ with $g(\rho+\mu)=0$ for all
$\mu\in\Lambda_+$ with $\ell(\mu)< d$. These spaces form a decreasing
filtration of $\cP^W$. We have
$$
E(g)(z)=\ell(z)g(z)-\sum_\eta f_\eta(z)g(z-\eta).
$$
Thus, \cite{CutOff} implies that each $U_d$ is $E$\_stable. Moreover,
$E-\ell(\rho)-d$ maps $U_d$ into $U_{d+1}$. This means that $E$ acts
on $U_d/U_{d+1}$ as scalar multiplication by $\ell(\rho)+d$. In
particular, $\ell(\rho)+d$ is not an eigenvalue of $E$ in
$U_{d+1}$. The action of $E$ is locally finite, since it preserves the
degree. This implies that there is a unique $E$\_stable complement
$\Uq_d$ of $U_{d+1}$ in $U_d$. Clearly, $\Uq_d$ is the (generalized)
eigenspace of $E$ in $\cP^W$ with eigenvalue $\ell(\rho)+d$. Since the
intersection of all $U_d$'s is $0$, there are no other other
eigenvalues. Thus $\cP^W=\oplus_d\Uq_d$.

Now assume that $\rho$ is not regular or dominant. For $N\in\NN$, let
$\cP^W_N$ be the $E$\_stable space of invariant polynomials of degree
$\le N$. For any $d$ choose an $N_d$ such that $\cP^W_N+U_d=\cP^W$ for
all $N\ge N_d$.

If $\rho$ is only dominant, then the map $\rho+\Lambda_+\pfeil V/W$ is
injective. Thus, the codimension of $U_d$ is independent of
$\rho$. This implies that for any $N\ge N_d$, the intersection
$U_d\cap\cP^W_N$ forms a family of subspaces of the finite dimensional
space $\cP^W_{N_d}$ which depends continuously on $\rho$. Also
$E$ depends continuously on $\rho$. It follows that
$U_d\cap\cP^W_N$, hence $U_d$ itself, is $E$\_stable. Then we conclude
as above.

Finally, for given $N$ and generic $\rho$ choose $d$ such that
$\cP^W_N\cap U_d=0$. Then $\cP^W_N$ is killed by $p(E)$, where
$p(z):=\prod_{i=0}^{d-1}(z-\ell(\rho)-i)$. Again by continuity,
$\cP^W_N$ is killed for all $\rho$. Since $p$ has no multiple zeros,
$E$ is diagonalizable on $\cP^W_N$, hence on $\cP^W$.\qed

\noindent So far, we did not use condition \co0.. Now, it will provide
the link between the function $\ell$ and the degree of a polynomial.

\Theorem Main-Prop-E. Let $\rho\in V_0$ be dominant.
\item{a)} For every $\lambda\in\Lambda_+$ there is a unique polynomial
$p_\lambda\in\cP^W$ with $\|deg|p_\lambda\le\ell(\lambda)$ and
$p_\lambda(\rho+\mu)=\delta_{\lambda\mu}$ (Kronecker delta) for all
$\mu\in\Lambda_+$ with $\ell(\mu)\le\ell(\lambda)$.
\item{b)}For every $d\in\NN$, the set of $p_\lambda$ with
$\ell(\lambda)\le d$ forms a basis of the space of $p\in\cP^W$ with
$\|deg|p\le d$.
\item{c)}The polynomial $p_\lambda$ is an eigenvector for $E$. More
precisely, $E(p_\lambda)=\ell(\rho+\lambda)p_\lambda$.
\Par

\Proof: We show first that {\it a)} implies {\it b)} and {\it c)}.

Let $\sum_\lambda a_\lambda p_\lambda=0$ be a non\_trivial linear
dependence relation. Choose $\lambda$ with $a_\lambda\ne0$ and
$\ell(\lambda)$ maximal. Then evaluation at $\rho+\lambda$ yields the
contradiction $a_\lambda=0$. Thus, the $p_\lambda$'s are linearly
independent.

Next, let $g\in\cP^W$ with $\|deg|g=d$. By induction we may assume
that {\it b)} holds for $d-1$. Hence there is a linear combination
$g'$ of $p_\mu$'s with $\ell(\mu)\le d-1$ such that $h:=g-g'$ vanishes
at all points $\rho+\mu$ with $\ell(\mu)<d$. Thus
$h':=h-\sum_{\ell(\lambda)=d} h(\rho+\lambda)p_\lambda$ vanishes in
all points $\rho+\mu$ with $\ell(\mu)\le d$. We have $\|deg|h'\le
d$. Thus, $h'\ne0$ contradicts the uniqueness of $p_\lambda$ with
$\ell(\lambda)=d$. Thus $g$ is a linear combination of the $p_\lambda$
with $\ell(\lambda)\le d$ which shows {\it b)}.

For {\it c)} it suffices, by continuity, to consider $\rho$ regular
and dominant. Consider the space $M_\lambda$ from \cite{PropL2}. By
{\it a)}, it contains, up to a scalar only one polynomial of degree
less at most $\ell(\lambda)$, namely $p_\lambda$. Thus, \cite{PropL2}
and \cite{PropE} imply that $p_\lambda$ is an eigenvector for $E$. A
direct calculation shows
$E(p_\lambda)(\rho+\lambda)=\ell(\rho+\lambda)$. Hence the eigenvalue
is $\ell(\rho+\lambda)$, which shows {\it c)}.

Now we prove {\it a)}. Condition \co0. guarantees the existence of one
regular dominant $\rho$ for which {\it a)}, hence {\it b)} and {\it
c)} hold. Let $\Lambda_+(d)$ be the set of $\mu\in\Lambda_+$ with
$\ell(\mu)\le d$. Let $\cP^W_d$ be the space of $g\in\cP^W$ with
$\|deg|g\le d$. Then {\it c)} shows in particular that
$\|dim|\cP^W_d=\#\Lambda_+(d)$. Thus $p_\lambda$ is defined by as many
(inhomogeneous) linear equations as there are variables. Its unique
solvability can be expressed by the non\_vanishing of a
determinant. This implies that {\it a)} holds for $\rho$ in the
complement of countably many hypersurfaces of $V_0$. This is, in
particular, a Zariski dense subset of $V_0$.

Now consider the action of $E$ on $\cP^W_d$. Then, by {\it c)},
$\prod_{i=0}^d(E-\ell(\rho)-i)$ is zero on $\cP^W_d$ for a Zariski
dense subset of $\rho$'s, hence for all $\rho\in V_0$. Let $F_i$ be
kernel of $E-\ell(\rho)-i$ in $\cP^W_d$. Then
$\cP^W_d=\oplus_{i=0}^dF_i$. The dimension of $F_i$ depends upper
semicontinuously on $\rho$. On the other hand their sum is
constant. Hence $\|dim|F_i$ is constant and equals the number $N_d$ of
$\mu\in\Lambda_+$ with $\ell(\mu)=d$.

Since $\rho$ is dominant, \cite{EigValE} implies that every
$g\in F(d)$ vanishes in $\rho+\Lambda_+(d-1)$. Moreover, the map
$F(d)\pfeil\CC^{N_d}: g\mapsto (g(\rho+\lambda)\mid\ell(\lambda)=d)$
is injective. Since both sides have the same dimension it is also
surjective. This implies that polynomials $p_\lambda$ as in {\it a)}
exist. Uniqueness follows again from the fact that the number of
equations equals the number of variables.\qed

\noindent We record this last fact for future reference:

\Corollary equaldim. Let $\rho\in V_0$ be dominant. For every $d$, the
dimension of the space of $g\in\cP^W$ with $\|deg|g\le d$ equals the
number of $\mu\in\Lambda_+$ with $\ell(\mu)\le d$.

\noindent
The equality $E(p_\lambda)=\ell(\rho+\lambda)p_\lambda$ can be
rewritten as
$$\eqno{MainE2}
\ell(z-\rho-\lambda)p_\lambda(z)=
\sum_{\eta\in\Lambda_1}f_\eta(z)p_\lambda(z-\eta).
$$
From this we obtain a formula for special values of $p_\lambda$. We
need:

\Definition: A {\it path\/} from $\lambda\in\Gamma$ to $\mu\in\Gamma$ is a
sequence $\tau_*=\tau_0,\tau_1,\ldots,\tau_d\in V$ with
$\tau_0=\lambda$, $\tau_d=\mu$ and $\tau_i-\tau_{i-1}\in\Lambda_1$ for
all $i=1,\ldots,d$. The path is {\it positive} if all $\tau_i$ are in
$\Lambda_+$.

\Definition: An element $\rho\in V_0$ is {\it non\_integral\/} if
$\alpha(\rho)\not\in\ZZ$ for all $\alpha\in\Delta$.

\medskip\noindent Observe that every $\rho$ coming from a
multiplicity free space is regular dominant but none is
non\_integral. Thus, certain formulas below are actually easier for
non\_geometric $\rho$\_shifts.

\Theorem Formula1. Let $\rho\in V_0$ be non\_integral. Let
$\lambda,\mu\in\Lambda_+$ with $d=\ell(\mu-\lambda)\ge0$. Then
$$\eqno{SpecialValue}
p_\lambda(\rho+\mu)=
{1\over d!}\sum_{\tau_*}
f_{\tau_1-\tau_0}(\rho+\tau_1)
f_{\tau_2-\tau_1}(\rho+\tau_2)\ldots
f_{\tau_d-\tau_{d-1}}(\rho+\tau_d),
$$
where the sum runs over all paths from $\lambda$ to
$\mu$. Moreover, only positive paths contribute to the sum. Thus, if
one restricts the sum to positive paths, then the formula is valid for
all regular dominant $\rho$.

\Proof: The non\_integrality of $\rho$ makes sure that none of the
denominators vanish. We proceed by induction on $d$. The statement
holds by definition for $d=0$. Let $d\ge1$. Putting $z=\rho+\mu$ in
\cite{MainE2} we obtain
$$
p_\lambda(\rho+\mu)=
{1\over d}\sum_{\eta\in\Lambda_1}f_\eta(\rho+\mu)
p_\lambda(\rho+\mu-\eta).
$$
By \cite{CutOff}, the coefficient $f_\eta(\rho+\mu)$ is zero whenever
$\tau_{d-1}=\mu-\eta$ is not in $\Lambda_+$. We conclude by
induction.\qed

\noindent As a corollary we get the extra vanishing property:

\Corollary ExtraVanishing. Let $\rho\in V_0$ be dominant. Let
$\Lambda$ be the monoid generated by $\Lambda_1$. Then
$p_\lambda(\rho+\mu)=0$ for every $\lambda,\mu\in\Lambda_+$ with
$\mu-\lambda\not\in\Lambda$.

\noindent For fixed $k\ge0$ we can sum over all paths with
$\tau_k=\tau$ first. Then we get

\Corollary. Assume $\rho\in V_0$ is regular dominant. Then for all
$\lambda,\mu\in\Lambda_+$ and all $k\in\NN$:
$$\eqno{ellPieri1}
{\ell(\mu-\lambda)\choose
k}p_\lambda(\rho+\mu)=\sum_{\tau\in\Lambda_+\atop\ell(\tau-\lambda)=k}
p_\lambda(\rho+\tau)p_\tau(\rho+\mu).
$$

\noindent
Observe that both sides of \cite{ellPieri1} depend polynomially on
$\mu$. Thus, if we set $\mu=z-\rho$, we obtain the following Pieri type
formula:

\Corollary.
$$\eqno{ellPieri2}
{\ell(z)-\ell(\rho+\lambda)\choose
k}p_\lambda(z)=\sum_{\tau\in\Lambda_+\atop\ell(\tau-\lambda)=k}
p_\lambda(\rho+\tau)p_\tau(z).
$$

\noindent
For $\lambda=0$ we get
\Corollary.
$$
{\ell(z)-\ell(\rho)\choose
k}=\sum_{\tau\in\Lambda_+:\ell(\tau)=k}p_\tau(z).
$$

\beginsection Pieri. Pieri rules

We consider the matrix of multiplication by $h\in\cP^W$ in the
$p_\lambda$\_basis:
$$\eqno{P-Formula 1}
h(z)p_\mu(z)=\sum_{\eta\in\Gamma}a_\eta^h(\mu)p_{\mu+\eta}(z)
$$
where we put $a_\eta^h(\mu)=0$ whenever $\mu+\eta\not\in\Lambda_+$.
We can compute the coefficients by evaluating both sides in the points
$z\in\rho+\Lambda_+$. The next proposition shows in particular that
the sum is over a finite set of $\eta$'s which is independent of $\mu$.

\Proposition Vanish2. $a_\eta^h(\mu)=0$ unless $\eta\in\Lambda$ and
$\ell(\eta)\le\|deg|h$.

\Proof: Fix $\mu\in\Lambda_+$ and choose $\eta_0\in\Gamma$ with
$\ell(\eta_0)$ minimal such that $\eta_0\not\in\Lambda$ but
$a_{\eta_0}^h(\mu)\ne0$. In particular,
$\mu+\eta_0\in\lambda_+$. Substituting $z=\rho+\mu+\eta_0$ in
\cite{P-Formula 1} we obtain by the extra vanishing property
(\cite{ExtraVanishing}) that
$$
0=h(\rho+\mu+\eta_0)p_\mu(\rho+\mu+\eta_0)=
a_{\eta_0}^h(\mu)+
\sum_{\eta\ne\eta_0}a_\eta^h(\mu)p_{\mu+\eta}(\rho+\mu+\eta_0).
$$
The $p$\_factor in the sum vanishes unless
$\eta_0-\eta\in\Lambda$. This implies $\eta\not\in\Lambda$ and
$\ell(\eta)<\ell(\eta_0)$. From this we get $a_\eta^h(\mu)=0$ by
minimality. Thus $a_{\eta_0}^h(\mu)=0$.

The inequality $\ell(\eta)\le\|deg|h$ simply reflects the fact that
the $p_\mu$ with $\ell(\mu)\le d$ form a basis of the space of
invariant polynomials of degree $\le d$ (\cite{Main-Prop-E}{\it b}).\qed

\Theorem Pieri1. Let $\rho\in V_0$ be non\_integral. Let
$\tau\in\Lambda$ with $d:=\ell(\tau)$ and $\mu\in\Lambda_+$
with $\mu+\tau\in\Lambda_+$. Then
$$\eqno{P-Formula 2}
a_\tau^h(\mu)=
\sum_{\tau_*}
\left[\sum_{i=0}^d{(-1)^{d-i}\over i!(d-i)!}h(\rho+\tau_i)\right]
f_{\tau_1-\tau_0}(\rho+\tau_1)\ldots f_{\tau_d-\tau_{d-1}}(\rho+\tau_d)
$$
where the outer sum runs over all paths from $\mu$ to
$\mu+\tau$. Moreover, only positive paths contribute to the sum. Thus, if
one restricts the sum to positive paths, then the formula is valid for
all regular dominant $\rho$.

\Proof: Substituting $z=\rho+\mu$ in \cite{P-Formula 1} yields
$a_0^h(\mu)=h(\rho+\mu)$ which is just \cite{P-Formula 2} for
$\tau=0$. Now we proceed by induction on $d$.

In \cite{P-Formula 1}, we substitute $z=\rho+\mu+\tau$ and obtain
$$
h(\rho+\mu+\tau)p_\mu(\rho+\mu+\tau)=
\sum_{\eta}a_\eta^h(\mu)p_{\mu+\eta}(\rho+\mu+\tau).
$$
The summands are zero unless $\eta,\tau-\eta\in\Lambda$. To all terms
with $\eta\ne\tau$ we can apply the induction hypothesis to the first
factor and \cite{SpecialValue} to the second. Thus we obtain (with
$\eta_i:=\tau_i-\tau_{i-1}$)
$$
a_\eta^h(\mu)p_{\mu+\eta}(\rho+\mu+\tau)=
$$
$$
=\sum_{\tau_*}\left[\sum_{i=0}^r{(-1)^{r-i}\over
(d-r)!i!(r-i)!}h(\rho+\tau_i)\right]
f_{\eta_1}(\rho+\tau_1)\ldots f_{\eta_d}(\rho+\tau_d)
$$
where $r=\ell(\eta)$ and the sum runs over all paths $\tau_*$ from
$\mu$ to $\mu+\tau$ with $\tau_r=\mu+\eta$. Thus we get
$$
\sum_{\eta\ne\tau}a_\eta^h(\mu)p_{\mu+\eta}(\rho+\mu+\tau)=
$$
$$
=\sum_{\tau_*}\left[\sum_{r=0}^{d-1}\sum_{i=0}^r{(-1)^{r-i}\over
(d-r)!i!(r-i)!}h(\rho+\tau_i)\right]
f_{\eta_1}(\rho+\tau_1)\ldots f_{\eta_d}(\rho+\tau_d),
$$
where we now sum over all paths from $\mu$ to $\mu+\tau$.
Now we interchange the order of summation in the bracket:
$$
\sum_{r=0}^{d-1}\sum_{i=0}^r{(-1)^{r-i}\over
(d-r)!i!(r-i)!}h(\rho+\tau_i)=
\sum_{i=0}^{d-1}\left[\sum_{r=i}^{d-1}{(-1)^{r-i}\over
(d-r)!i!(r-i)!}\right]h(\rho+\tau_i).
$$
The sum in brackets can be rewritten as
$$
\sum_{r=i}^{d-1}{(-1)^{r-i}\over(d-r)!i!(r-i)!}=
{1\over(d-i)!i!}\sum_{r=i}^{d-1}(-1)^{r-i}{d-i\choose r-i}=
-{(-1)^{d-i}\over (d-i)!i!}
$$
Assembling everything together we get
$$
a_\tau^h(\mu)=h(\rho+\tau_d)p_\mu(\rho+\tau_d){+}\sum_{\tau_*}
\left[\sum_{i=0}^{d-1}{(-1)^{d-i}\over i!(d{-}i)!}h(\rho{+}\tau_i)\right]
f_{\eta_1}(\rho{+}\tau_1)\ldots f_{\eta_d}(\rho{+}\tau_d).
$$
By \cite{SpecialValue}, the first summand is nothing but the
missing  case $i=d$ of the second one. This yields \cite{P-Formula 2}.\qed

\Corollary. Let $\rho\in V_0$ be dominant. Let $\tau\in\Lambda$
and $\mu\in\Lambda_+$ with $\mu+\tau\in\Lambda_+$. Then
$$
a_\tau^h(\mu)=\sum_\eta(-1)^{\ell(\tau-\eta)}h(\rho+\mu+\eta)
p_\mu(\rho+\mu+\eta)p_{\mu+\eta}(\rho+\mu+\tau)
$$
where the sum runs over all $\eta\in\Lambda$ with
$\tau-\eta\in\Lambda$ and $\mu+\eta\in\Lambda_+$.

\Proof: By continuity, we may assume that $\rho$ is non\_integral. In
\cite{P-Formula 2}, we fix an $i$ and sum over all different
$\eta:=\tau_i$ first. Then the formula follows from two applications
of \cite{Formula1}.\qed

\beginsection Operators. Construction of other difference operators

We are going to need the following

\Lemma Phi+. Let $\tau\in\Lambda$ and $\omega\in\Phi$ with
$\omega(\tau)<0$. Then $-\omega\in\Phi$.

\Proof: Since $\tau\in\Lambda$ there is $\eta\in\Lambda_1$ with
$\omega(\eta)<0$. By definition, there is a $w\in W$ with
$w\eta\in\Sigma$. Then $w^{-1}\omega\not\in\Phi^+$ which implies that
$-w^{-1}\omega$, hence $-\omega$ is in $\Phi$.\qed

Now we can prove a cut-off property which is dual to that in
\cite{CutOff}.

\Lemma DualCutOff. Let $\rho\in V_0$ be regular dominant,
$\mu\in\Lambda_+$, $\eta\in\Lambda_1$, and $\lambda:=\mu+\eta$. Then
$f_\eta(-\rho-\mu)=0$ whenever $\lambda\not\in\Lambda_+$.

\Proof: Since $\lambda\not\in\Lambda_+$ there is $\omega\in\Sigma$
with $\omega(\lambda)=\omega(\mu)+\omega(\eta)<0$. Hence,
\co4. implies $\omega(\mu)=0$ and $\omega(\eta)=-1$. Since \cite{Phi+}
implies $\overline\omega:=-\omega\in\Phi$, the factor
$\overline\omega(z)-k_{\overline\omega}=-\omega(z)-\omega(\rho)$ of
$f_\eta(z)$ vanishes at $z=-\rho-\mu$.\qed

For $d\in\NN$ we define the falling factorial functions as
$$
[z\da d]:=z(z-1)\ldots(z-d+1).
$$
Now we generalize the definition of $f_\tau$ to all $\tau\in\Lambda$ as
follows:
$$\eqno{Def-f-l}
f_\tau(z):={
\prod\limits_{\omega\in\Phi:\omega(\tau)>0}[\omega(z)-k_\omega\da\omega(\tau)]
\over
\prod\limits_{\alpha\in\Delta:\alpha(\tau)>0}[\alpha(z)\da\alpha(\tau)]}.
$$
\noindent Now we need stronger non\_degeneracy conditions for $\rho$.

\Definition: An element $\rho\in V_0$ is called {\it strongly
dominant\/} if $\rho$ is regular dominant and
$\omega(\rho)-k_\omega\not\in\ZZ_{<0}$ and
$\omega(\rho)+k_\omega\not\in\ZZ_{\le0}$ for all $\omega\in\Phi^+$.

\medskip\noindent
The set of strongly dominant $\rho$ is Zariski\_dense in $V_0$. In
fact, it suffices to show that it is non\_empty since then it is the
complement of countably many hyperplanes. To produce a strongly
dominant $\rho$, we set all $k_\omega$ equal $t$ where
$t\not\in\QQ_{\le0}$. Every $\alpha\in\Delta^+$ and $\omega\in\Phi^+$
has an expression $\sum_{\omega_i\in\Sigma}a_i\omega_i$ with all
$a_i\in\ZZ_{\ge0}$. Let $N:=\sum_ia_i\in\ZZ_{>0}$. Then
$\alpha(\rho)=Nt\not\in\ZZ_{\le0}$,
$\omega(\rho)+k_\omega=(N+1)t\not\in\ZZ_{\le0}$, or
$\omega(\rho)-k_\omega=(N-1)t\not\in\ZZ_{<0}$.

In fact, one can show that for the examples coming from multiplicity
free spaces $\rho$ is strongly dominant whenever all $k_\omega$ are
real and positive.

\Lemma nichtVerschwind. Let $\rho\in V_0$ be strongly dominant
and $\lambda,\mu\in\Lambda_+$, $\tau:=\lambda-\mu$. Then
$f_\tau(\rho+\lambda)$ and $f_\tau(-\rho-\mu)$ are defined and
non\_zero.

\Proof: Let $\alpha\in\Delta$ with $\alpha(\tau)>0$. Consider the
factor
$F=[\alpha(\rho+\lambda)\da\alpha(\tau)]=
\big(\alpha(\rho)+\alpha(\lambda)\big)\ldots
\big(\alpha(\rho)+\alpha(\mu)+1\big)$.
If $\alpha\in\Delta^+$, then $\alpha(\lambda)\ge\alpha(\mu)+1>0$. From
$\alpha(\rho)\not\in\ZZ{\le0}$ we get $F\ne0$. If $-\alpha\in
\Delta^+$, then $0\ge\alpha(\lambda)\ge\alpha(\mu)+1$ and
$\alpha(\rho)\not\in\ZZ_{\ge0}$ which again implies
$F\ne0$. Therefore, $f_\tau(\rho+\lambda)$ is defined.

The denominator of $f_\tau(-\rho-\mu)$ as well as the numerators are
treated analogously by
considering the factors
$$
[\alpha(-\rho-\mu)\da\alpha(\tau)]=
\pm\big(\alpha(\rho)+\alpha(\lambda)-1\big)
\ldots\big(\alpha(\rho)+\alpha(\mu)\big),
$$
$$
[\omega(\rho+\lambda)-k_\omega\da\omega(\tau)]=
\big(\omega(\rho)-k_\omega+\omega(\lambda)\big)\ldots
\big(\omega(\rho)-k_\omega+\omega(\mu)+1\big),
$$
$$
[\omega(-\rho-\mu)-k_\omega\da\omega(\tau)]=\pm
\big(\omega(\rho)+k_\omega+\omega(\lambda)-1\big)\ldots
\big(\omega(\rho)+k_\omega+\omega(\mu)\big).
$$\qed

\noindent In particular, for every $\lambda\in\Lambda_+$ we can define the
{\it virtual dimension\/} as
$$
d_\lambda:=(-1)^{\ell(\lambda)}
{f_\lambda(-\rho)\over f_\lambda(\rho+\lambda)}\ne0.
$$
The terminology comes from the fact that, for multiplicity free
spaces, $d_\lambda$ is indeed the dimension of the simple module with
highest weight $\lambda$. This will be proved in a forthcoming
paper. In general, the virtual dimension can be rewritten as
$$
d_\lambda=\prod_{\alpha\in\Delta^+}{\alpha(\rho+\lambda)\over\alpha(\rho)}\ 
\prod_{\omega\in\Phi^+}{[\omega(\rho+\lambda)+k_\omega-1\da\omega(\lambda)]
\over[\omega(\rho+\lambda)-k_\omega\da\omega(\lambda)]}.
$$
The first factor is just Weyl's dimension
formula. It is easy to see that $d_\lambda$ is a polynomial function
in $\lambda$ if and only if $k_\omega\in{1\over2}\ZZ_{>0}$ for all
$\omega$. In this case we have
$$
d_\lambda=\prod_{\alpha\in\Delta^+}{\alpha(\rho+\lambda)\over\alpha(\rho)}\ 
\prod_{\omega\in\Phi^+}\prod_{s=-k_\omega+1}^{k_\omega-1}
{\omega(\rho+\lambda)+s\over\omega(\rho)+s}.
$$
Another property of the virtual dimension is:

\Theorem Dual1. Let $\rho$ be strongly dominant. Let
$\lambda,\mu\in\Lambda_+$ with $\tau=\lambda-\mu\in\Lambda$. Then
$$\eqno{FormulaDual}
{d_\lambda\over d_\mu}=(-1)^{\ell(\tau)}
{f_\tau(-\rho-\mu)\over f_\tau(\rho+\lambda)}.
$$

\Proof: With
$$
A_\alpha:={
[-\alpha(\rho+\mu)\da\alpha(\tau)]
\over
[\alpha(\rho+\lambda)\da\alpha(\tau)]},
\quad
B_\omega:={[-\omega(\rho+\mu)-k_\omega\da\omega(\tau)]
\over
[\omega(\rho+\lambda)-k_\omega\da\omega(\tau)]},
$$
we get
$$\eqno{FormulaAB}
{f_\tau(-\rho-\mu)\over f_\tau(\rho+\lambda)}=
\prod_{\alpha\in\Delta:\alpha(\tau)>0}A_\alpha^{-1}
\prod_{\omega\in\Phi:\omega(\tau)>0}B_\omega.
$$
It is customary to extend the definition of the falling factorial
functions as $[z\da d]:=1/(z+1)\ldots(z-d)$ when $d<0$. With this
convention, the formula $[z\da d]=1/[z-d\da -d]$ holds for all
$d\in\ZZ$ if the left-hand side is non-zero. This holds in our case and
we obtain
$$
{[-\omega(\rho+\mu)-k_\omega\da\omega(\tau)]
\over
[\omega(\rho+\lambda)-k_\omega\da\omega(\tau)]}
=
{[\omega(\rho+\mu)-k_\omega\da-\omega(\tau)]
\over
[-\omega(\rho+\lambda)-k_\omega\da-\omega(\tau)]}.
$$
Now I claim that
$$
\Phi_{\tau>0}=\Phi^+_{\tau>0}\cup(-\Phi^+_{\tau<0})
$$
where the subscript $\tau>0$ means ``the subset of all $\omega$
with $\omega(\tau)>0$'', etc. If $\omega\in\Phi$, then either
$\omega\in\Phi^+$ or $-\omega\in\Phi^+$ (by \co 2.) which shows the
inclusion ``$\subseteq$''. Conversely, let $\omega\in\Phi^+$ with
$\omega(\eta)<0$. Then $-\omega\in\Phi_{\tau>0}$ by \cite{Phi+},  which
proves the claim.

We have $k_\omega=k_{-\omega}$ and therefore $B_\omega=B_{-\omega}$
whenever both $\pm\omega\in\Phi$. Thus, the claim implies
that in \cite{FormulaAB} we can replace the product over all
$\omega\in\Phi_{\tau>0}$ by the product over all
$\omega\in\Phi^+$. Similarly, we can change the range of the first
product to $\Delta^+$.

Now we apply the formulas
$$
[z\da a-b]={[z+b\da a]\over[z+b\da b]}={[z\da a]\over[z-(a-b)\da b]}
$$
with $a=\omega(\lambda)$ and $b=\omega(\mu)$. They hold whenever none
of the denominators are zero. That this is so in our case follows from
\cite{nichtVerschwind}. We obtain
$$
B_\omega=
{[-\omega(\rho)-k_\omega\da\omega(\lambda)]
\over
[-\omega(\rho)-k_\omega\da\omega(\mu)]}
{[\omega(\rho+\mu)-k_\omega\da\omega(\mu)]
\over
[\omega(\rho+\lambda)-k_\omega\da\omega(\lambda)]}.
$$
Similarly, we have
$$
A_\alpha={
[-\alpha(\rho)\da\alpha(\mu)]
\over
[-\alpha(\rho)\da\alpha(\lambda)]
}{
[\alpha(\rho+\lambda)\da\alpha(\lambda)]
\over
[\alpha(\rho+\mu)\da\alpha(\mu)]
}
$$
for all $\alpha\in\Delta^+$. Since $\lambda,\mu\in\Lambda_+$ we can
multiply in the definition \cite{Def-f-l} of $f_\lambda$, $f_\mu$ over all
$\alpha\in\Delta^+$ and $\omega\in\Phi^+$. The asserted formula
\cite{FormulaDual} follows readily.\qed

For $\tau\in\Lambda$ with $d:=\ell(\tau)$ we now define the rational function
$$\eqno{Def-b-t}
b_\tau^h(z):=\sum_{\tau_*}
\left[\sum_{i=0}^d{(-1)^{d-i}\over i!(d-i)!}h(z-\tau_i)\right]
f_{\tau_1-\tau_0}(z-\tau_0)\ldots f_{\tau_d-\tau_{d-1}}(z-\tau_{d-1})
$$
where the sum runs over all paths from $0$ to
$\tau$. For $h\in\cP$ define $h^-\in\cP$ by $h^-(z)=h(-z)$.

\Theorem Vanish. Let $\rho$ be strongly dominant and
non\_integral. Let $\lambda,\mu\in\Lambda_+$ with
$\tau=\lambda-\mu\in\Lambda$. Then
$$\eqno{PEquation}
b_\tau^h(-\rho-\mu)=(-1)^{\ell(\tau)}
{d_\lambda\over d_\mu} a_\tau^{h^-}(\mu).
$$

\Proof: This would follow directly from \cite{Pieri1} and \cite{Dual1}
if in \cite{Def-b-t} the summation
were restricted to those paths $\tau_*$ for which $\mu+\tau_*$ is
positive.  But other paths do not contribute to $b_\tau^h(-\rho-\mu)$
anyway: let $i$ be minimal such that
$\mu+\tau_i\not\in\Lambda_+$. Then
$f_{\eta_i}(-\rho-\mu-\tau_{i-1})=0$, by \cite{DualCutOff}.\qed

\noindent The following consequence is crucial:

\Corollary Crucial. Let $\tau\in\Lambda$ with $\ell(\tau)>\|deg|h$. Then
$b_\tau^h(z)=0$.

\Proof: For strongly dominant, non\_integral $\rho$ this follows from
\cite{Vanish} and \cite{Vanish2} since the set of points $z=-\rho-\mu$
with $\mu\in\Lambda_+$ and $\tau+\mu\in\Lambda_+$ is Zariski
dense. For general $\rho$ we conclude by continuity.\qed

\noindent Thus, for each $h\in\cP^W$ we can define the difference operator
$$
D_h:=\sum_{\tau\in\Lambda}b_\tau^h(z)T_\tau.
$$
We are going to rewrite it: Fix $\tau\in\Lambda$ and put
$i:=\ell(\tau)$. Then (with $\eta_i:=\tau_i-\tau_{i-1}$)
$$
\sum_{\tau_*,\tau_i=\tau}h(z-\tau_i)f_{\eta_1}(z-\tau_0)\ldots
f_{\eta_d}(z-\tau_{d-1})=
$$
$$
=\sum_{\tau_*,\tau_i=\tau}(T_{\eta_1}\ldots T_{\eta_i})(h)\cdot
f_{\eta_1}\cdot
T_{\eta_1}(f_{\eta_2})\cdot
(T_{\eta_1}T_{\eta_2})(f_{\eta_3})\cdot\ldots\cdot
(T_{\eta_1}\ldots T_{\eta_{d-1}})(f_{\eta_d}).
$$
This is easily recognized as the coefficient of $T_\tau$ in
$L^ihL^{d-i}$ where we regard $h$ as  a multiplication operator. Thus we
get
$$
\sum_{\ell(\tau)=d}b_\tau^h(z)T_\tau
={1\over d!}\sum_{i=0}^d(-1)^{d-i}{d\choose i}L^ihL^{d-i}=
{1\over d!}(\|ad|L)^d(h).
$$
Using \cite{Crucial} we get

\Corollary. Let $h\in\cP^W$ and $d\in\NN$ with $d>\|deg|h$. Then
$(\|ad|L)^d(h)=0$.

\noindent We also obtained the first part of:

\Theorem Main. a) Let $h\in\cP^W$. Then $D_h=\|exp|(\|ad|L)(h)$.\Par
\noindent
b) The map $h\mapsto D_h$ is an algebra homomorphism. In particular,
the $D_h$ commute pairwise.\Par\noindent
c) If $\rho$ is dominant, then the action of the $D_h$ on $\cP^W$ is
simultaneously diagonalizable. More precisely, if
$\lambda\in\Lambda_+$ then $D_h(p_\lambda)=h(\rho+\lambda)p_\lambda$.

\Proof: Let
$\cA\subseteq\|End|(\cP^W)$ be the largest subalgebra on which
$\|ad|L$ acts locally nilpotently. Then $\cP^W\subseteq\cA$. Moreover,
$\|ad|L$ is a derivation, hence $\|exp|(\|ad|L)$ is an automorphism of
$\cA$. This shows {\it b)}.

In {\it c)}, we may assume that $\rho$ is regular dominant and then
conclude by continuity. First, observe $E=D_\ell$. Therefore $D_h$
and $E$ commute. The space of $f\in\cP^W$ with $\|deg|f\le e$ can be
characterized as the direct sum of the $E$\_eigenspaces for the
eigenvalues $\ell(\rho), \ell(\rho)+1,\ldots,\ell(\rho)+e$. Therefore,
$D_h(f)\le\|deg|f$ for all $f\in\cP^W$. On the other hand, both $L$
and $h$, hence $D_h$, preserve the space $M_\lambda$ from
\cite{PropL2}. We conclude that $D_h(p_\lambda)=cp_\lambda$ for some
constant $c$. The constant term of $D_h$, i.e., the coefficient of
$T_0$, is $h$. Thus, evaluating at $z=\rho+\lambda$ gives
$c=h(\rho+\lambda)$. This shows c).\qed

\beginsection Further. Further analysis of the difference operators

In this section, we derive some basic properties of the functions
$b_\tau^h(z)$. First the degree:

\Proposition b-degree. Let $h\in\cP^W$ and $\tau\in\Lambda$. Then
$\|deg|b_\tau^h(z)\le\|deg|h$.

\Proof: Let $c(z)$ be a rational function of degree $d$ on $V$ and
$\tau\in\Gamma$. Then $[L,cT_\tau]=[L,c]T_\tau+c[L,T_\tau]$. We have
$[L,c]=\sum_\eta f_\eta(z)(c(z-\tau)-c(z))$. Since
$\|deg|f_\eta(z)\le1$ and $\|deg|(c(z-\tau)-c(z))<d$ we see that
$[L,c]T_\tau$ has coefficients of degree $\le d$. Similarly for the
other term: The coefficients in  $[L,T_\tau]=\sum_\eta
c(z)(f_\eta(z)-f_\eta(z))T_{\tau+\eta}$ have degree $\le d$. Thus we
have shown that $\|ad|L$ does not increase the degrees of
coefficients. The assertion follows from the formula \cite{Main}a).\qed

\noindent Next we study the denominator:

\Proposition b-denom. Let $h\in\cP^W$. For fixed
$\tau\in\Lambda$ with $b_\tau^h\ne0$ and $\alpha\in\Delta$
put
$$
S(\alpha,\tau):=\{i\in\ZZ\mid i\ne0, b_{\tau+i\alpha^\vee}^h\ne0\}.
$$
Then
$$
b_\tau^h(z)\prod_{\alpha\in\Delta^+}
\prod_{i\in S(\alpha,\tau)}
(\alpha(z-\tau)-i)
$$
is a polynomial in $z$.

\Proof: By the explicit formula in \cite{Main}{\it a)}, the only
denominators which can occur are products of terms $\alpha(z)-m$ with
$\alpha\in\Delta$ and $m\in\ZZ$. Fix one such factor $\alpha(z)-m$ and
let $S$ be the set of all $\tau$ such that $\alpha(z)-m$ occurs in the
denominator of $b_\tau$. Let $\cH_0$ be the hyperplane $\{\alpha-m=0\}$.

Fix $\tau\in S$ and let $z\in\cH_0$. Suppose that none of the points
$z-\tau'$, $\tau\in S$, $\tau'\ne\tau$ is in the $W$\_orbit of
$z-\tau$. Then one could find a {\it symmetric\/} function $f$ which
does not vanish in $z-\tau$ but vanishes to a very high order in all
other points $z-\tau'$. Then $D_h(f)$ would not be regular in $z$.
Thus, there must be $\tau'$ and $w$ with $w(z-\tau')=z-\tau$. Since
there are only finitely many choices of $w$ and $\tau'$ there is one
choice which works for a dense subset of points $z\in\cH_0$. By
continuity, we get
$$\eqno{UUU}
w(z-\tau')=z-\tau\hbox{ for all }z\in\cH_0.
$$
Choose any $z_0\in\cH_0$. Then for any $y$ with $\alpha(y)=0$ and any
number $t$ we have $z=t*y+z_0\in\cH_0$. Comparing the coefficient of
$t$ in \cite{UUU} yields $w(y)=y$. Since $w$ cannot be the identity we
get $w=s_\alpha$. Also $\tau'$ is unique since
$$
\tau'=z-s_\alpha(z-\tau)=\tau+(m-\alpha(\tau))\alpha^\vee.
$$
This yields $\alpha(z)-m=\alpha(z-\tau)-i$ with $i=m-\alpha(\tau)\in
S(\alpha,\tau)$.

It remains to show that $\alpha(z)-m$ occurs with multiplicity one in
the denominator of $b_\tau^h$. Let $N>0$ be the larger of the powers
with which $\epsilon:=\alpha(z)-m$ occurs in the denominator of
$b_\tau^h$ and $b_{\tau'}^h$. Put $c:=\epsilon^Nb_\tau^h$ and likewise
$c':=\epsilon^Nb_{\tau'}^h$. Then for every $f\in\cP^W$, the rational
function $c(z)f(z-\tau)+c'(z)f(z-\tau')$ has a zero of order at least
$N$ along the divisor $\epsilon(z)=0$.  From the identity
$z-\tau'=s_\alpha(z-\tau-\epsilon\alpha^\vee)$ we infer that
$$\eqno{VVV}
c(z)f(z-\tau)+c'(z)f(z-\tau')=
$$
$$
\qquad=(c(z)+c'(z))f(z-\tau)+c'(z)(f(z-\tau-\epsilon\alpha^\vee)-f(z-\tau)).
$$
We conclude that $c(z)+c'(z)$ is divisible by $\epsilon$. Hence, since
one of $c$ or $c'$ is not divisible by $\epsilon$ the other is not
either. Next observe that $\tau\ne\tau'$ implies that $\cH_0-\tau$ is
not the reflection hyperplane of $s_\alpha$. Hence there exists $z\in
cH_0$ and $f\in\cP^W$ which is divisible by $\epsilon$ such that the
directional derivative $D_{\alpha^\vee}f(z-\tau)\ne0$. This
implies that the right-hand side of \cite{VVV} vanishes to exactly
first order in $\cH_0$. Thus $N=1$.\qed

\noindent A lower bound for the numerator is given by

\Proposition b-num. Let $\rho$ be non\_integral. Then the numerator of
$b_\tau^h$ is divisible by
$$
\prod_{\omega\in\Phi:\omega(\tau)>0}[\omega(z)-k_\omega\da\omega(\tau)].
$$

\Proof: The non\_integrality of $\rho$ ensures that the denominator of
$b_\tau^h(z)$ does not vanish whenever $z\in\rho+\Gamma$. If
$\omega\in\Sigma$, then the definition of $b_\tau^h(z)$ along with
\cite{DualCutOff} implies that $b_\tau^h(z)=0$ for all
$z=\rho+\lambda$ with $\lambda\in\Lambda_+$ and
$\omega(\lambda-\tau)<0$. Thus, $b_\tau^h$ is divisible by
$[\omega(z)-k_\omega\da\omega(\tau)]$.

Let $\omega\in\Phi$ be arbitrary with $\omega(\tau)>0$. There is $w\in
W$ with $w\omega\in\Sigma$. Thus, by the case above,
$b_{w\tau}^h(z)$ is divisible by $[w\omega(z)-k_\omega\da
\omega(\tau)]$. The rest follows from the fact that $D_h$ is a
symmetric operator: $b_{w\tau}^h(wz)=b_\tau(z)$.\qed

We now introduce an order relation on $\Lambda$: We define
$\tau_1\le\tau_2$ if $\ell(\tau_1)<\ell(\tau_2)$ or
$\ell(\tau_1)=\ell(\tau_2)$ and $\tau_2-\tau_1$ is a sum of positive
roots.

\Theorem leading coeff. For $h\in\cP^W$ let $\tau\in\Lambda_+$ be
maximal (with respect to ``$\le$'') with $b_\tau^h\ne0$. Then
$\ell(\tau)=\|deg|h$, $\tau\in\Lambda_+$ and
$b_\tau^h(z)\in\CC^*f_\tau(z)$.

\Proof: First of all, it follows from \cite{Vanish2} and \cite{Vanish}
that $\ell(\tau)=\|deg|h$. The set of $\tau$ with $b_\tau^h\ne0$ is
$W$\_stable. It follows that maximality of $\tau$ implies
$\alpha(\tau)\ge0$ for all $\alpha\in\Delta^+$. From
$s_\alpha(\tau)=\tau-\alpha(\tau)\alpha^\vee$ it follows that
$i=-\alpha(\tau)$ is the minimal element of $S(\alpha,\tau)$ and all
other $i$ satisfy $-\alpha(\tau)\le i<0$. Thus Propositions~\ncite{b-denom} and
\ncite{b-num} imply that $c(z):=b_\tau^h(z)f_\tau(z)^{-1}$ is a
polynomial. Moreover, by \cite{b-degree} we have
$\|deg|c(z)\le\|deg|h-\|deg|f_\tau=\|deg|\ell(\tau)-\|deg|f_\tau=-j(\tau)$
where $j(\tau):=\|deg|f_\tau-\ell(\tau)$. Now all assertions follow
from the following claim: {\it Let $\tau\in\Lambda$ be dominant for
$\Delta^+$. Then $j(\tau)\ge0$ and equality holds if and only if
$\tau\in\Lambda_+$.}

To prove the claim let $a_+:=\|max|(a,0)$ for any $a\in\ZZ$. Then,
from the definition of $f_\tau$ it follows that
$$
\|deg|f_\tau=\sum_{\omega\in\Phi}\omega(\tau)_+
-\sum_{\alpha\in\Delta^+}\alpha(\tau).
$$ Thus $\tau\mapsto j(\tau)$ is a piecewise linear convex
function. Moreover, $j(\tau)=0$ whenever $\tau\in\Lambda_+$ (by
\co2. and \co3.). Together this implies $j(\tau)\ge0$ for all $\tau$.
Every $\omega_0\in\Sigma$ defines a codimension\_1\_face of
$\Lambda_+$. Let $\tau\in\Lambda$ be close enough to that face such
that $\omega(\tau)\ge0$ for all $\omega\in\Phi^+$ except for
$\omega=\omega_0$ where $\omega_0(\tau)<0$. Since $-\omega_0\in\Phi$
by \cite{Phi+}, \co3. then implies
$j(\tau)=0-\omega_0(\tau)+(-\omega_0)(\tau)=-2\omega_0(\tau)>0$.
Therefore, $j(\tau)$ takes strictly positive values outside
$\Lambda_+$ which finishes the proof of the claim.\qed

\Corollary highterm. Assume $\rho$ is strongly dominant. For
$\lambda\in\Lambda_+$ let $h=p_\lambda^-$. Then $b_\tau^h=0$ unless
$\tau\le\lambda$ and $b_\lambda^h=f_\lambda(\rho+\lambda)^{-1}f_\lambda$.

\Proof: Let $\tau$ be maximal with $b_\tau^h\ne0$. Then
$\tau\in\Lambda_+$ and $b_\tau^h=cf_\tau$ for some $c\ne0$. For
$g:=h^-=p_\lambda$ we get by definition
$a_\tau^g(0)=\delta_{\lambda\tau}$. \cite{Vanish} implies
$$
cf_\tau(-\rho)=b_\tau^h(-\rho)=
(-1)^{\ell(\tau)}d_\tau\delta_{\lambda\tau}=
{f_\tau(-\rho)\over f_\tau(\rho+\tau)}\delta_{\lambda\tau}.
$$
Hence \cite{nichtVerschwind} implies
$\tau=\lambda$ and $c=f_\lambda(\rho+\lambda)^{-1}$.\qed

With our Pieri formula we can convert this into a triangularity
result. For this, it is convenient to introduce another normalization
of the $p_\lambda(z)$: let
$P_\lambda:=f_\lambda(\rho+\lambda)p_\lambda$. Thus
$P_\lambda(\lambda+\rho)=f_\lambda(\rho+\lambda)$.

\Corollary. For all $\lambda,\mu\in\Lambda_+$ holds
$$
P_\lambda P_\mu=
P_{\lambda+\mu}+\sum_{\nu<\lambda+\mu}c_{\lambda\mu}^\nu P_\nu.
$$

\Proof: For $h=P_\lambda^-$ we have just seen $b_\tau^h=0$ unless
$\tau\le\lambda$ and $b_\lambda^h=f_\lambda$. \cite{Vanish2}
implies $P_\lambda P_\mu=\sum_\nu c_{\lambda\mu}^\nu P_\nu$ with
$c_{\lambda\mu}^\nu=f_\lambda(\rho+\lambda){f_\mu(-\rho)\over
f_\nu(-\rho)}b_{\nu-\mu}^h(-\rho-\mu)$. Thus $c_{\lambda\mu}^\nu=0$ unless
$\nu-\mu\le\lambda$. Moreover,
$$
c_{\lambda\mu}^{\lambda+\mu}=
{f_\lambda(-\rho-\mu)f_\mu(-\rho)\over f_{\lambda+\mu}(-\rho)}.
$$
This last expression equals 1 as follows from the identity $[z-b\da
a][z\da b]=[z\da a+b]$.\qed

\noindent In the classical and semiclassical case, $P_\lambda$
is exactly the polynomial obtained by normalizing the leading
coefficient to 1. To make sense of this in general we introduce a
monomial basis of $\cP^W$ as follows: consider
$\Sigma=\{\omega_1,\ldots,\omega_r\}$ and its dual basis
$\Sigma^\vee=\{\eta_1,\ldots,\eta_r\}$. Then for any
$\lambda\in\Lambda_+$ we define
$$
{\bf e}_\lambda:=\prod_{i=1}^r P_{\eta_i}^{\omega_i(\lambda)}.
$$
Then we get easily by induction:

\Corollary. For every $\lambda\in\Lambda_+$ there is an expansion
$$
P_\lambda={\bf e}_\lambda+\sum_{\mu<\lambda}d_{\lambda\mu}{\bf e}_\mu.
$$

\beginsection multiplicity free. Multiplicity free spaces

In this section, we introduce the main class of examples to which the
theory developed in the preceding sections applies.

Let $G$ be a connected reductive group (everything is defined over
$\CC$) and $U$ a finite dimensional $G$\_module. Let $\cO(U)$ be its
algebra of polynomial functions. Then $U$ is called a {\it
multiplicity free space\/} if $\cO(U)$ is a multiplicity free
$G$\_module, i.e., every simple $G$\_module occurs in $\cO(U)$ at most
once. A more geometric criterion is due to Vinberg\_Kimelfeld
(\cite{ViKi}, see also \cite{Montreal} Thm. 3.1): a Borel subgroup of
$G$ has a dense orbit in $U$.

We assume from now on that $U$ is a multiplicity free space and we are going
to derive a structure $(\Gamma,\Sigma,W,\ell)$ from it. For a dominant
integral weight $\lambda$ let $M^\lambda$ be the simple $G$\_module
with {\it lowest\/} weight $-\lambda$. Let $\Lambda_+$ be the set
of $\lambda$ such that $M^\lambda$ occurs in $\cO(U)$. Thus, as a
$G$\_module, we have
$$
\cO(U)\cong\oplus_{\lambda\in\Lambda_+} M^\lambda.
$$
We regard characters as elements of the dual Cartan algebra
$\ft^*$. Let $\Gamma^\vee\subseteq\ft^*$ be the subgroup
generated by $\Lambda_+$. Then we can define the first ingredient as
$\Gamma:=\|Hom|(\Gamma^\vee,\ZZ)$.

It is known (\cite{HoUm}, see also \cite{Montreal} Thm. 3.2) that
$\Gamma_+$ is a monoid which is generated by a basis $\Sigma^\vee$ of
$\Gamma^\vee$. Let $\Sigma\subseteq\Gamma$ be its dual basis. Since
$U$ is a vector space, the algebra $\cO(U)$ is graded and every
irreducible constituent $M^\lambda$ occurs in some degree
$\ell(\lambda)$. It is easy to see that $\ell$ is additive on
$\Lambda_+$. Hence it extends to a linear function $\ell\in\Gamma$.

\def\rhoq{{\overline\rho}}

The reflection group $W$ is more involved to construct. Let $\cD(U)$
be the algebra of polynomial coefficient differential operators on
$U$. We are interested in the algebra of $G$\_invariant operators
$\cD(U)^G$. Every $D\in\cD(U)^G$ acts on $M^\lambda\subseteq\cO(U)$ by
a scalar, denoted by $c_D(\lambda)$. Let $V\subseteq\ft^*$ be the
$\CC$\_span of $\Gamma^\vee$. By \cite{Montreal} Cor.~4.4, the
function $c_D(\lambda)$ extends to a polynomial function $c_D(z)$ on
$V$. Thus we get an injective homomorphism
$c:\cD(U)^G\into\cO(V):D\mapsto c_D$.

To determine its image let $\rhoq\in\ft^*$ be the half\_sum of the
positive roots and $\Wq\subseteq GL(\ft^*)$ the Weyl group of $G$. The
twisted action of $\Wq$ on $\ft^*$ is defined as
$w\bullet\chi:=w(\chi+\rhoq)-\rhoq$. Then $W$ is characterized as
follows:

\Theorem Harish Chandra. There is a unique subgroup $W\subseteq\Wq$ such
that
\item{a)}the subspace $V$ and the lattice $\Gamma^\vee$ are stable
under the twisted action of $W$;
\item{b)}the image of $c$ consists exactly of the invariants under
this twisted $W$\_action.\Par

\noindent
This finishes the description of the structure
$(\Gamma,\Sigma,W,\ell)$. The main point is the following theorem
whose proof will occupy the rest of this section.

\Theorem MFAxiom. Let $(\Gamma,\Sigma,W,\ell)$ be the structure
derived from a multiplicity free space. Then all axioms
\co11. through \co0. hold.

{\everypar{\leavevmode\hangindent20pt}

\noindent \co11. follows from the fact that $\cD(U)^G$ is a polynomial
ring (\cite{HoUm}; see also \cite{Montreal} Cor.~4.7) and the
Shepherd-Todd theorem.

\noindent
\co1. is clear since all weights in $\Sigma^\vee$ are dominant.

\noindent
\co3''. follows from \cite{Harish Chandra} since the Euler vector
field $\xi$ is in $\cD(U)^G$ and we have $c_\xi=\ell$.

\noindent
\co3'. is trivial, since degrees of non\_constant polynomials are
strictly positive.

\noindent
\co5. Any linear $W$\_invariant $f$ comes from a $G$\_invariant
differential operator of order one, hence from a $G$\_invariant vector
field $\xi$ on $U$. We have $U^\vee\subseteq\cO(U)$ and
$\xi(U^\vee)\subseteq U^\vee$. Thus, $\xi$ is uniquely determined by
the $G$\_endomorphism $\xi_0:=\xi|_U$. This endomorphism $\xi_0$ acts
on each simple component of lowest weight $\eta$ of $U$ as scalar
$f(\eta)$. But these lowest weights run exactly through
$\Sigma_1^\vee$.

}

\noindent \co2., \co3., and \co4. are handled case by case. For this,
we first need some reductions. Assume that there is a reductive group
$\Gq$ with $G\subseteq\Gq\subseteq GL(U)$ such that $G$ is normal in
$\Gq$ and the quotient $\Gq/G$ is a torus. Then the center $\Zq$ of
$\Gq$ acts as a scalar on each module $M^\lambda\subseteq\cO(U)$. This
means that $U$ is also multiplicity free with respect to $\Gq$ and
that there is an isomorphism $\overline\Lambda_+=\Lambda_+$. A
differential operator is in $\cD(U)^G$ if and only if it acts as a
scalar on each $M^\lambda$. This shows that
$\cD(U)^G=\cD(U)^\Gq$. Hence we obtain an isomorphism
$(\Gamma,\Sigma,W,\ell)\cong
(\overline\Gamma,\overline\Sigma,\Wq,\overline\ell)$.

This observation is applied as follows: let $U=U_1\oplus\ldots\oplus
U_s$ be the decomposition of $U$ into simple modules and let
$A=\G_m^s\subseteq GL(U)$ consisting of the scalar multiplications in
each factor. Then we can replace $G$ by $\Gq=AG$, i.e., assume right
away that $A\subseteq G$. A multiplicity free space with that property is
called {\it saturated}.

If $(G_1,U_1)$ and $(G_2,U_2)$ are multiplicity free spaces then
$(G,U)=(G_1\times G_2, U_1\oplus U_2)$ is one as well. If $U_1$ and
$U_2$ are non\_zero, then $U$ is called {\it decomposable}. The
combinatorial structures are related by
$(\Gamma,\Sigma,W,\ell)=
(\Gamma_1\oplus\Gamma_2,\Sigma_1\cup\Sigma_2,W_1\times W_2,\ell_1+\ell_2)$.
Moreover, one readily verifies that all axioms, but in particular
\co2., \co3., and \co4., hold for $U$ whenever they hold for $U_1$ and
$U_2$.

Thus we may assume that $U$ is indecomposable and saturated. These
multiplicity free spaces have been classified independently by
Benson\_Ratcliff, \cite{BenRat1}, and Leahy, \cite{Leahy}. This
classification together with the structure $(\Gamma,\Sigma,W,\ell)$ is
tabulated in the next section from which one easily verifies the three
axioms case by case.

Finally, we verify axiom \co0.. This will also provide the motivation
for the whole theory. Observe, that, with $U$, the dual
representation $U^\vee$ is also a multiplicity free space. Its algebra of
functions decomposes as a $G$\_module like
$$
\cO(U^\vee)=\oplus_{\lambda\in\Lambda_+}M_\lambda
$$
where $M_\lambda=(M^\lambda)^\vee$ is the simple $G$\_module with
{\it highest\/} weight $\lambda$. An element $D\in\cO(U^\vee)$ can be
regarded as a differential operator with constant coefficients on
$U$. Let $\cD(U)$ be the algebra of all polynomial coefficient
differential operators. Thus we get a $G$\_module isomorphism
$$
\cO(U)\otimes\cO(U^\vee)\pfeil\cD(U):f\otimes D\mapsto fD.
$$
Using this isomorphism we can construct a distinguished basis of
$\cD(U)^G$ as follows: the space of $G$\_fixed vectors in
$M^\lambda\otimes M_\mu$ is non\_zero if and only if $\lambda=\mu$ and
in that case it is one\_dimensional. Let $D_\lambda\in\cD(U)^G$ be the
image of a generator. Then the family $D_\lambda$,
$\lambda\in\Lambda_+$, is a basis of $\cD(U)^G$. Thus, the polynomials
$c_\lambda(z):=c_{D_\lambda}(z)$ form a basis of the space of shifted
invariant polynomials on $V$. To get rid of the shift, choose any
$W$\_equivariant projection $\pi:\ft^*\pfeil V$ and any vector
$\sigma\in V^W$. Let $\kappa:=\rhoq-\pi(\rhoq)$ and
$\rho:=\pi(\rhoq)+\sigma=\rhoq-\kappa+\sigma\in V$. Then the fact that
$\rhoq+V$ is $W$\_stable means precisely that $\kappa$ is
$W$\_fixed. Hence the shifted $W$\_action with $\rhoq$\_shift
coincides with that corresponding to the $\rho$\_shift. Thus, if we
put $p_\lambda(z):=c_\lambda(z-\rho)$ then we obtain a basis of the
(unshifted) $W$\_invariants on $V$. Now, condition \co0. is
essentially Theorem~4.10 of \cite{Montreal}. That theorem also makes sure
that we can normalize $D_\lambda$ in such a way that
$p_\lambda(\rho+\lambda)=c_\lambda(\lambda)=1$.

The only point left to show is that $\pi$ and $\sigma$ can be chosen
in such a way that $\rho\in V_0$. We are even going to construct a
canonical element $\rho\in V_0$.

Recall the following consequence of the local structure theorem (see,
e.g., \cite{Montreal} Thm.~2.4). Let $B^-$ be a Borel subgroup
opposite to $B$. Then there is a point $u\in U$ such that $B^-u$ is
open in $U$. Moreover, there is a parabolic subgroup $P\subseteq G$
with Levi part $L$ and unipotent radical $R_uP$ such that
\item{$\bullet$}the orbit $Pu=B^-u$;

\item{$\bullet$}the isotropy group $P_u$ is contained and normal in
$L$;

\item{$\bullet$}the quotient $L/P_u$ is a torus.

\noindent
For $\lambda\in\Lambda_+$, the lowest weight vector of $f_\lambda\in
M^\lambda$ can be normalized to $f_\lambda(u)=1$ and then defines a
character of $P$, hence of $L$. The intersection of the kernels of
these characters equals $L_u=P_u$. Therefore, we can identify $V$ with
the dual of $\|Lie|L/P_u$. Furthermore, we can choose a subspace $C$
of the center of $\|Lie|L$ such that $\|Lie|L=C\oplus\|Lie|P_u$ and
$V\cong C^*$. For a Cartan subalgebra $\ft$ of $\|Lie|L$ we obtain
$\ft^*=C^*\oplus\ft_u^*$. Now we choose for $\pi$ the projection of
$\ft^*$ onto $C^*=V$.

Let $w_L$ be the longest element of the Weyl group of $L$ with respect
to $\ft$. Then clearly $\pi(\rhoq)=\pi(w_L\rhoq)$, i.e.,
$\tilde\rho:={1\over2}(\rhoq+w_L\rhoq)$ has the same image in $V$ as
$\rhoq$. Let $\chi\in\ft^*$ be sum of all weights of $U$. Since it
comes from a character of $\|Lie| G$, it is $\Wq$\_invariant. Let
$\rho:={1\over2}\chi+\tilde\rho={1\over2}(\chi+\rhoq+w_L\rhoq)$. Then
$\pi(\rho)$ differs from $\pi(\rhoq)$ by the $W$\_invariant element
$\sigma:={1\over2}\pi(\chi)$.

Now I claim that $\pi(\rho)=\rho$, i.e., $\rho$ is already in
$V$. Consider the action of $\ft_u$ on the top exterior power of the
tangent space $\Lambda^{\|top|}T_uU=\Lambda^{\|top|}U$. On one side,
this is just the restriction of $\chi$ to $\ft_u$. On the other side,
observe that, $-2\tilde\rho$ is the sum of all roots belonging to
$R_uP$. Moreover, $T_uU=\|Lie|Pu=\|Lie|R_uP\oplus C$ with trivial
action of $\ft_u$ on $C$. Thus, the restriction of $\rho$ to $\ft_u$
is zero which means $\rho\in V$.

This element $\rho$ is very easy to calculate in every given case. For
the indecomposable saturated multiplicity free spaces it is recorded
in the tables below. This shows in particular that $\rho\in V_0$ in
every given case and concludes the proof of \cite{MFAxiom}.

An immediate consequence is the following statement. It would be
desirable to have a conceptual proof.

\Corollary. The polynomials $p_\lambda$ describing the spectrum of Capelli
operators on a multiplicity free space are the joint eigenfunctions of
a family of commuting difference operators.

\beginsection Tables. Tables

The table below lists all structures $(\Gamma,\Sigma,W,\ell)$ which
come from indecomposable saturated multiplicity free actions. First,
the combinatorial structure is defined (indicated by a double line
$||$) and then its relation to multiplicity free spaces.

The space $V$ is a subspace of some $\CC^m$ with canonical basis $e_i$
and coordinates $z_i$. In all cases, except {\bf III}$_{\rm odd}$ and
{\bf IVa}, we have $V=\CC^m$. In case {\bf V}, we found it easier to
work with basis vectors $e_i$, $e_i'$, $e''$ and corresponding
coordinates $z_i$, $z_i'$, $z''$.

The Weyl group is given as follows: $s_{ij}$ denotes the transposition
$z_i\leftrightarrow z_j$. The notation $S_3(z_1,z_3,z_5)$ means the
symmetric group permuting the coordinates $z_1,z_3,z_5$ and leaving
all others fixed. A similar convention holds for other reflection
groups, e.g., ${\Sans D}_3(z_2,z_4,z_6)$. Finally, $\pm z_i$ means the
reflection about the hyperplane $z_i=0$.

The sets $\Delta^+$ and $\Phi^+$ are given such that condition
\co3. can be verified easily.

We also give the orbit structure of $\Sigma$ under the group $\pm
W$. Each entry corresponds to an element of $\Sigma$. Then
$\omega_i\in\pm W\omega_j$ if and only if the $i$-th and $j$-th entry
are equal, disregarding the sign. The sign ``$\pm$'' stands if and
only if $-\omega_i\in W\omega_i$. Otherwise, the same or
different sign means that $\omega_i\in W\omega_j$ or $-\omega_i\in
W\omega_j$, respectively.

Then we give the $z$\_coordinates of a general element $\rho\in
V_0$. If in the $\pm W$-table $\omega\in\Sigma$ has an upper\_case
letter $R$, $S$, etc. then $\omega(\rho)$ is denoted by the
corresponding lower case letter $r$, $s$, etc.

Below the definition of the structure we list the indecomposable
saturated multiplicity free actions giving rise to it. The list is
comprehensive by the classification in \cite{Kac}, \cite{BenRat1}, and
\cite{Leahy}. In \cite{Montreal} we collected all the necessary data to
verify the assertions in the table.

In each case, we list the values of $\omega(\rho)$, $\omega\in\Sigma$
where $\rho=\1(\chi+\rhoq+w_L\rhoq)$ is the canonical choice of $\rho$
as described in the preceding section. Then we describe how the basis
vectors $e_i$ are related to actual weights. The notation is quite
straightforward: $\epsilon_i$ denotes a weight in the defining
representation of a classical group, $\alpha_i$ and $\omega_i$ are
simple roots and fundamental weights (numbered as in Bourbaki
\cite{Bou}). Weights of different factors of $G$ are distinguished by
primes: $\epsilon_i$, $\epsilon_i'$, $\epsilon_i''$, etc.

\medskip \def\\{\hfill\break} \def\1{{\textstyle{1\over2}}}
\def\phan#1{\setbox1=\hbox{#1}\wd1=0pt\box1}

\long\def\lstrich#1{\setbox0=\vbox{#1}
\penalty10000\medskip\penalty10000\noindent
\phan{\vrule height \ht0\hskip.5mm\vrule height \ht0}
\box0
\bigskip}

\noindent{\bf Case I: The classical cases:} ($1\le n$) (see also \cite{SymCap})

\lstrich{
$\Sigma:=\{z_1-z_2,z_2-z_3,\ldots,z_{n-1}-z_n,z_n\}$.

$\Sigma^\vee=\{e_1,e_1+e_2,e_1+e_2+e_3,\ldots,e_1+e_2+\ldots+e_n\}$

$\ell:=z_1+z_2+\ldots+z_n$

$W:=S_n=\<s_{12},s_{23},\ldots,s_{n{-1}\,n}\>$

$\Delta^+=\{z_i-z_j\mid1\le i<j\le n\}$

$\Phi^+=\{z_i-z_j\mid1\le i<j\le n\}\cup
\{z_i\mid1\le i\le n\}$

$\pm W$-orbits of $\Sigma$: $[\pm R,\pm R,\ldots,\pm R,S]$

$\rho=((n-1)r+s,(n-2)r+s,\ldots,r+s,s)$}

\noindent
$GL_p(\CC)$ on $S^2(\CC^p)$ with $1\le p$

$n=p$, $r=\1$, $s=\1$, $e_i=2\epsilon_i$

\medskip\goodbreak
\noindent
$GL_p(\CC)\times GL_q(\CC)$ on $\CC^p\otimes\CC^q$ with $1\le p\le q$

$n=p$, $r=1$, $s=\1(q-p+1)$, $e_i=\epsilon_i+\epsilon_i'$

\medskip\goodbreak
\noindent
$GL_p(\CC)$ on $\Lambda^2(\CC^p)$ with $2\le p$

$\vcenter{
\halign{#\ \hfill&#\ \hfill&#\ \hfill&#\hfill,\ 
$e_i=\epsilon_{2i-1}+\epsilon_{2i}$\cr
$p$ even:&$n={p\over2}$,  &$r=2$,&$s=\1$\cr
$p$ odd: &$n={p-1\over2}$,&$r=2$,&$s={3\over2}$\cr}}
$

\medskip\goodbreak
\noindent
$Sp_{2p}(\CC)$ on $\CC^{2p}$ with $1\le p$

$n=1$, $r$ undefined, $s=p$, $e_1=\epsilon_1$

\medskip\goodbreak
\noindent
$SO_p(\CC)\times\CC^*$ on $\CC^p$ with $3\le p$

$n=2$, $r={p\over2}-1$, $s=\1$,
$e_1=\epsilon_1+\epsilon$, $e_2=-\epsilon_1+\epsilon$

\goodbreak
\noindent
$Spin_{10}(\CC)\times\CC^*$ on $\CC^{16}$

$n=2$, $r=3$, $s={5\over2}$,
$e_1=\omega_5+\epsilon$, $e_1+e_2=\omega_1+2\epsilon$

\medskip\goodbreak
\noindent
$Spin_7(\CC)\times\CC^*$ on $\CC^8$

$n=2$, $r=3$, $s=\1$,
$e_1=\omega_3+\epsilon$, $e_2=-\omega_3+\epsilon$

\medskip\goodbreak
\noindent
${\rm G}_2\times\CC^*$ on $\CC^7$

$n=2$, $r={5\over2}$, $s=\1$,
$e_1=\omega_1+\epsilon$, $e_2=-\omega_1+\epsilon$

\medskip\goodbreak
\noindent
${\rm E_6}\times\CC^*$ on $\CC^{27}$

$n=3$, $r=4$, $s=\1$,
$e_1=\omega_1+\epsilon$, $e_1+e_2=\omega_6+2\epsilon$, $e_1+e_2+e_3=3\epsilon$

\medskip\hrule\goodbreak\medskip
\noindent{\bf Case II: The semiclassical cases:} ($3\le n$) (see also
\cite{Semisym})
\lstrich{
$\Sigma:=\{z_1-z_2,z_2-z_3,\ldots,z_{n-1}-z_n,z_n\}$

$\Sigma^\vee=\{e_1,e_1+e_2,e_1+e_2+e_3,\ldots,e_1+e_2+\ldots+e_n\}$

$\ell:=z_1+z_3+z_5+\ldots$

$W:=\{\pi\in S_n\mid\forall i:\pi(i)-i\ \|even|\}=
\<s_{13},s_{24},\ldots,s_{n{-2}\,n}\>$

$\Delta^+=\{z_i-z_j\mid1\le i<j\le n, i-j \ \|even|\}$

$\Phi^+=\{z_i-z_j\mid1\le i<j\le n, i-j \ \|odd|\}
\cup\{z_i\mid1\le i\le n, n-i\ \|even|\}$

$\pm W$-orbits of $\Sigma$: $[R,-R,R,-R,\ldots,S]$

$\rho=((n-1)r+s,(n-2)r+s,\ldots,r+s,s)$}

\noindent
$GL_p(\CC)\times GL_q(\CC)$ on $(\CC^p\otimes\CC^q)\oplus\CC^q$ with
$1\le p$, $2\le q$

$p<q$: $n=2p+1$, $r=\1$, $s={q-p\over2}$, $e_{2i}=\epsilon_i$ ($i=1,\ldots,p$),
$e_{2i-1}=\epsilon_i'$ ($i=1,\ldots,p+1$)

$p\ge q$: $n=2q$, $r=\1$, $s={p-q+1\over2}$,
$e_{2i}=\epsilon_i$ ($i=1,\ldots,q$),
$e_{2i-1}=\epsilon_i'$ ($i=1,\ldots,q$)
\medskip\goodbreak

\noindent
$GL_1(\CC)\times GL_q(\CC)$ on $(\CC\otimes\CC^q)\oplus(\CC^q)^*$ with
$2\le q$

$n=3$, $r={q-1\over2}$, $s={1\over2}$, $e_1=-\epsilon'_n$,
$e_2=\epsilon+\epsilon_1'+\epsilon_n'$, $e_3=-\epsilon_1'$.

\medskip\goodbreak

\noindent
$GL_p(\CC)$ on $\Lambda^2(\CC^n)\oplus\CC^n$ with $3\le p$

$n=p$, $r=1$, $s=\1$, $e_i=\epsilon_i$



\medskip\hrule\goodbreak\medskip
\noindent{\bf Case III: The quasiclassical cases:} ($3\le n$)

\lstrich{
\noindent\ \ {\bf $n$ odd}

$V:=\{z\in\CC^{n+1}\mid z_n=0\}$

$\Sigma:=\{z_1-z_2,z_2-z_3,\ldots,z_n-z_{n+1}\}$

$\Sigma^\vee=\{e_1,e_1+e_2,\ldots,e_1+e_2+\ldots+e_{n-1},-e_{n+1}\}$

$\ell:=\sum_{i=1}^{n+1}\1(1-3(-1)^i)z_i=2z_1-z_2+2z_3-+\ldots-z_{n+1}$

$W:=\{\pi\in S_{n+1}\mid\forall i:\pi(i)-i\ \|even|,\pi(n)=n\}=
\<s_{13},s_{24},\ldots,s_{n{-3}\,n{-}1},s_{n{-1}\,{n+1}}\>$

$\Delta^+=\{z_i-z_j\mid1\le i<j\le n+1, i-j \ \|even|,j\ne n\}$

$\Phi^+=\{z_i-z_j\mid1\le i<j\le n+1, i-j \ \|odd|\}$
(with $z_n=0$)

$\pm W$-orbits of $\Sigma$: $[R,-R,R,-R,\ldots,R,S,-S]$

$\rho=((n-2)r+s,(n-3)r+s,\ldots,r+s,s,0,-s)$}
\medskip
\lstrich{
\noindent\ \ {\bf $n$ even}

$\Sigma:=\{z_1-z_2,z_2-z_3,\ldots,z_{n-1}-z_n,z_{n-1}\}$

$\Sigma^\vee=\{e_1,e_1+e_2,\ldots,e_1+e_2+\ldots+e_{n-2},-e_n,
e_1+e_2+\ldots+e_n\}$

$\ell:=\sum_{i=1}^n\1(1-3(-1)^i)z_i=2z_1-z_2+2z_3-+\ldots-z_n$

$W:=\{\pi\in S_n\mid\forall i:\pi(i)-i\ \|even|\}=
\<s_{13},s_{24},\ldots,s_{n{-2}\,n}\>$

$\Delta^+=\{z_i-z_j\mid1\le i<j\le n, i-j \ \|even|\}$

$\Phi^+=\{z_i-z_j\mid1\le i<j\le n, i-j \ \|odd|\}
\cup\{z_i\mid 1\le i\le n, i\ \|odd|\}$

$\pm W$-orbits of $\Sigma$: $[R,-R,R,-R,\ldots,R,S]$

$\rho=((n-2)r+s,(n-3)r+s,\ldots,r+s,s,-r+s)$}

\noindent {\it Remark:} For $n=3$ and $n=4$, the structures in {\bf
II} and {\bf III} are isomorphic.

\medskip\noindent
$GL_p(\CC)\times GL_q(\CC)$ on $(\CC^p\otimes\CC^q)\oplus(\CC^q)^*$ with
$1\le p$, $2\le q$

$p<q$: $n=2p+1$, $r=\1$, $s={q-p\over2}$, $e_{2i{-}1}{=}\epsilon_i$
($i{=}1,\ldots,p$), $e_{2i}{=}\epsilon_i'$ ($i{=}1,\ldots,p$),
$e_{n+1}{=}\epsilon_q'$

$p\ge q$: $n=2q$, $r=\1$, $s={p-q+1\over2}$,
$e_{2i-1}=\epsilon_i$ ($i=1,\ldots,q$),
$e_{2i}=\epsilon_i'$ ($i=1,\ldots,q$)
\medskip\goodbreak
\noindent
$GL_p(\CC)\times\CC^*$ on $\Lambda^2(\CC^p)\oplus(\CC^p)^*$ with $3\le p$

$n=p$, $r=1$, $s=\1$, $e_i=\epsilon_i$ ($i=1,\ldots,2\lfloor{n\over2}\rfloor$),
$e_{n+1}=\epsilon_n$ if $n$ is odd


\medskip\hrule\goodbreak\medskip\goodbreak
\noindent{\bf Case IVa:}
\lstrich{
$V:=\{z\in\CC^7\mid z_2+z_4+z_6+z_7=0\}$

$\Sigma:=\{z_1-z_2,z_2-z_3,z_3-z_4,z_4-z_5,z_5-z_6,z_4+z_6\}$

\leavevmode\vbox{\halign{$#$\hfill&$#$\hfill\cr
\vrule height \baselineskip width 0pt
\Sigma^\vee=\{&(1,0,0,0,0,0,0),(1,1,0,0,0,0,-1),(1,1,1,0,0,0,-1),\cr
&(\1,\1,\1,\1,-\1,-\1,-\1),(\1,\1,\1,\1,\1,-\1,-\1),(\1,\1,\1,\1,\1,\1,
-{3\over2})\}\cr}}

$\ell:=2(z_1+z_3+z_5)$

$W:=S_3(z_1,z_3,z_5)
\times S_4(z_2,z_4,z_6,z_7)$

$\Delta^+=\{z_i-z_j\mid (i,j)=(1,3),(1,5),(3,5),
(2,4),(2,6),(2,7),(4,6),(4,7),(6,7)\}$

\leavevmode\vbox{\halign{$#$\hfill&$#$\hfill\cr
\vrule height \baselineskip width 0pt
\Phi^+=&\{\|sign|(j-i)(z_i-z_j)\mid i=1,3,5; j=2,4,6,7\}\cup\cr
&\{(z_i+z_j)\mid (i,j)=(2,4),(2,6),(4,6)\}\cr}}

$\pm W$-orbits of $\Sigma$: $[R,-R,R,-R,R,\pm S]$

$\rho=({s\over2}+4r,{s\over2}+3r,{s\over2}+2r,{s\over2}+r,{s\over2},
{s\over2}-r,-{3\over2}s-3r)$}

\noindent
$Sp_{2p}(\CC)\times GL_3(\CC)$ on $\CC^{2p}\otimes\CC^3$ with $3\le p$

$r=\1$, $s=p-2$,
$e_1=\epsilon_1'+\epsilon_2'$, $e_3=\epsilon_1'+\epsilon_3'$,
$e_5=\epsilon_2'+\epsilon_3'$,\\\indent
$e_2-e_7=\epsilon_1+\epsilon_2$,
$e_4-e_7=\epsilon_1+\epsilon_3$,
$e_6-e_7=\epsilon_2+\epsilon_3$


\medskip\hrule\goodbreak\medskip
\noindent{\bf Case IVb:}

\lstrich{
$\Sigma:=\{\1(z_1-z_2-z_4-z_6),z_2-z_3,z_3-z_4,z_4-z_5,z_5-z_6,z_5+z_6\}$

\leavevmode\vbox{\halign{$#$\hfill&$#$\hfill\cr
\vrule height \baselineskip width 0pt
\Sigma^\vee=\{&(2,0,0,0,0,0),(1,1,0,0,0,0),(1,1,1,0,0,0),(2,1,1,1,0,0),\cr
&({3\over2},\1,\1,\1,\1,\1),(\1,\1,\1,\1,\1,-\1)\}\cr}}

$\ell:=2z_1$

$W:={\Ss D}_3(z_2,z_4,z_6)
\times{\Ss C}_2(z_3,z_5)$

$\Delta^+=\{z_i\pm z_j\mid (i,j)=(2,4),(2,6),(4,6)\}
\cup\{z_3\pm z_5,2z_3,2z_5\}$

\leavevmode\vbox{\halign{$#$\hfill&$#$\hfill\cr
\vrule height \baselineskip width 0pt
\Phi^+=&\textstyle
\{z_i\pm z_j\mid2\le i<j\le6, i-j\hbox{ odd}\}
\cup\cr
&\textstyle\{\1(z_1\pm z_2\pm z_4\pm z_6)\mid
\hbox{1 or 3 minus signs}\}\cr}}

$\pm W$-orbits of $\Sigma$: $[S,\pm R,\pm R,\pm R,\pm R,\pm R]$

$\rho=(2s+6r,4r,3r,2r,r,0)$}

\noindent
$Sp_4(\CC)\times GL_p(\CC)$ on $\CC^4\otimes\CC^p$ with $4\le p$

$r=\1$, $s={p-3\over2}$, $e_3=\epsilon_1+\epsilon_2,
e_5=\epsilon_1-\epsilon_2$

$2e_1{=}\epsilon_1'{+}\epsilon_2'{+}\epsilon_3'{+}\epsilon_4',
2e_2{=}\epsilon_1'{+}\epsilon_2'{-}\epsilon_3'{-}\epsilon_4',
2e_4{=}\epsilon_1'{-}\epsilon_2'{+}\epsilon_3'{-}\epsilon_4',
2e_6{=}{-}\epsilon_1'{+}\epsilon_2'{+}\epsilon_3'{-}\epsilon_4'$

\medskip\hrule\goodbreak\medskip
\noindent{\bf Case IVc:}

\lstrich{
$\Sigma:=\{z_1-z_2,z_2-z_3,z_3-z_4,z_4-z_5,z_4+z_5\}$

$\Sigma^\vee=\<(1,0,0,0,0),(1,1,0,0,0),(1,1,1,0,0),
               (\1,\1,\1,\1,\1),(\1,\1,\1,\1,-\1)\>$

$\ell:=2(z_1+z_3+z_5)$

$W:=S_3(z_1,z_3,z_5)
\times{\Ss C}_2(z_2,z_4)$

$\Delta^+=\{z_i-z_j\mid
(i,j)=(1,3),(1,5),(3,5)\}\cup\{z_2\pm z_4,2z_2,2z_4\}$

$\Phi^+=
\{z_i\pm z_j\mid1\le i<j\le5, i-j \hbox{ odd}\}$

$\pm W$-orbits of $\Sigma$: $[R,-R,R,-R,R]$

$\rho=(4r,3r,2r,r,0)$}

\noindent
$Sp_4(\CC)\times GL_3(\CC)$ on $\CC^4\otimes\CC^3$

$r=\1$, $e_2=\epsilon_1+\epsilon_2,
e_4=\epsilon_1-\epsilon_2$, $e_1=\epsilon_1'+\epsilon_2',
e_3=\epsilon_1'+\epsilon_3',
e_5=\epsilon_2'+\epsilon_3'$


\medskip\hrule\goodbreak\medskip
\noindent{\bf Case V:} ($1\le b\le a\le3$)

\lstrich{%
\leavevmode\vbox{\halign{$#$\hfill&$#$\hfill\cr
\Sigma=\{&z_1-z_2,\ldots,z_{a-1}-z_a,z'_1-z'_2,\ldots,z'_{b-1}-z'_b,\cr
&z_a+z'_b-z'',z_a-z'_b+z'',-z_a+z'_b+z''\}\cr}}

\leavevmode\vbox{\halign{$#$\hfill&$#$\hfill\cr
\vrule height \baselineskip width 0pt
\Sigma^\vee=\{&e_1,e_1+e_2,\ldots,e_1+e_2+\ldots+e_{a-1},
                e_1',e_1'+e_2',\ldots,e_1'+e_2'+\ldots+e_{b-1}',\cr
&\1\sum e_i+\1\sum e_i',\1\sum e_i+\1e'',\1\sum e_i'+\1e''\}\cr}}

$\ell:=2(z_1+z'_1)$

$W:=(\ZZ/2\ZZ)^{a+b-1}=\{(z_1,\pm z_2,\ldots,\pm z_a,
z'_1,\pm z'_2,\ldots,\pm z'_b,\pm z'')\}$

$\Delta^+=\{2z_2,\ldots,2z_a,2z'_2,\ldots,2z'_b,2z''\}$

\leavevmode\vbox{\halign{$#$\hfill&$#$\hfill\cr
\vrule height \baselineskip width 0pt
\Phi^+=&\{z_i\pm z_{i+1}\mid1\le i<a\}\cup
\{z'_i\pm z'_{i+1}\mid1\le i<b\}\cup\cr
&\{\pm z_a\pm z'_b\pm z''\mid\hbox{at most one minus sign}\}\cr}}

\leavevmode\vbox{\halign{#\hfill&#\hfill\cr
\vrule height \baselineskip width 0pt
$\pm W$-orbits of $\Sigma$: &$a\ge b>1$:
$[R_1,\pm R_2,\ldots,\pm R_{a-1},R_1',\pm R_2',\ldots,\pm R_{b-1}',
\pm S,\pm S,\pm S]$\cr
&$a>b=1$: $[R_1,\pm R_2,\ldots,\pm R_{a-1},S,-S,S]$\cr}}

$\rho=(r_1+\ldots+r_{a-1}+s,r_2+\ldots+r_{a-1}+s,
\ldots,r_{a-1}+s,s$;\\\indent
\phantom{$\rho=($}
$r'_1+\ldots+r'_{b-1}+s,r'_2+\ldots+r'_{b-1}+s,\ldots,
r'_{b-1}+s,s;s)$ (for $a>b\ge1$)}

\noindent
{\it Remark:} The cases $(a,b)=(1,1)$ and $(2,1)$ are isomorphic to
the cases $n=3$ and $n=4$ of {\bf Case II}, respectively.
\medskip

\noindent
$(Sp_{2p}(\CC)\times\CC^*)\times GL_2(\CC)$
on $(\CC^{2p}\otimes\CC^2)\oplus\CC^2$ with $2\le p$

$a=3$, $b=1$, $r_1=\1$, $r_2=p-1$, $s=\1$
\\\indent
$e_1=2\epsilon+\omega_2'$, $e_2=\omega_2$, $e_3=\alpha_1$, $e_1'=\omega_2'$,
$e''=\alpha_1'$

\medskip\goodbreak
\noindent
$GL_p(\CC)\times SL_2(\CC)\times GL_q(\CC)$ on
$(\CC^p\otimes\CC^2)\oplus(\CC^2\otimes\CC^q)$ with $2\le p,q$

$a=2$, $b=2$, $r_1={p-1\over2}$, $r_1'={q-1\over2}$, $s=\1$
\\\indent
$e_1=\omega_2$, $e_2=\alpha_1$, $e_1'=\omega_2''$, $e_2'=\alpha_1''$,
$e''=2\omega'$

\medskip\goodbreak
\noindent
$(Sp_{2p}(\CC)\times\CC^*)\times SL_2(\CC)\times GL_q(\CC)$ on
$(\CC^{2p}\otimes\CC^2)\oplus(\CC^2\otimes\CC^q)$ with $2\le p,q$

$a=3$, $b=2$, $r_1=\1$, $r_2=p-1$, $r_1'={q-1\over2}$, $s=\1$
\\\indent
$e_1=2\epsilon$, $e_2=\omega_2$, $e_3=\alpha_1$, $e_1'=\omega_2''$,
$e_2'=\alpha_1''$, $e''=2\omega'$
\medskip\goodbreak
\noindent
$(Sp_{2p}(\CC)\times\CC^*)\times SL_2(\CC)\times (Sp_{2q}(\CC)\times\CC^*)$ on
$(\CC^{2p}\otimes\CC^2)\oplus(\CC^2\otimes\CC^{2q})$ with $2\le p,q$

$a=3$, $b=3$, $r_1=\1$, $r_2=p-1$, $r_1'=\1$, $r_2'=q-1$, $s=\1$
\\\indent
$e_1=2\epsilon$, $e_2=\omega_2$, $e_3=\alpha_1$, $e_1'=2\epsilon''$,
$e_2'=\omega_2''$, $e_3'=\alpha_1''$, $e''=2\omega'$

\medskip\goodbreak
\noindent
$Spin_8(\CC)\times\CC^*\times\CC^*$ on $\CC^8_+\oplus\CC^8_-$

$a=2$, $b=2$, $r_1=\1$, $r_1'=\1$, $s={3\over2}$
\\\indent
$e_1=2\epsilon$, $e_2=\epsilon_1-\epsilon_4$,$e_1'=2\epsilon'$,
$e_2'=\epsilon_1+\epsilon_4$, $e''=\epsilon_2+\epsilon_3$

\medskip\hrule\goodbreak\medskip

\noindent{\bf Case VIa:}

\lstrich{
$\Sigma:=\{z_1-z_2,z_2-z_3,2z_3\}$

$\Sigma^\vee=\{(1,0,0),(1,1,0),(\1,\1,\1)\}$

$\ell:=2z_1$

$W:=(\ZZ/2\ZZ)^2=\{(z_1,z_2,z_3)
\mapsto(z_1,\pm z_2,\pm z_3)\}$

$\Delta^+=\{2z_2,2z_3\}$

$\Phi^+=\{z_1\pm z_2,z_2\pm z_3,2z_3\}$

$\pm W$-orbits of $\Sigma$: $[R,\pm S,\pm T]$

$\rho=(r+s+{t\over2},s+{t\over2},{t\over2})$}

\noindent
$Sp_{2p}(\CC)\times GL_2(\CC)$ on $\CC^{2p}\otimes\CC^2$ with $2\le p$

$r=\1$, $s=p-1$, $t=1$,
$e_1=\omega_2'$, $e_2=\omega_2$, $e_3=\alpha_1+\alpha_1'$

\medskip\goodbreak
\noindent
$Spin_9(\CC)\times\CC^*$ on $\CC^{16}$

$r=\1$, $s=2$, $t=3$,
$e_1=2\epsilon$, $e_2=\omega_1$, $e_3=-\omega_1+2\omega_4$

\medskip\hrule\goodbreak\medskip
\noindent{\bf Case VIb:}

\lstrich{
$\Sigma:=\{z_1-z_2,z_2-z_3,z_3-z_4,z_3+z_4\}$

$\Sigma^\vee=\{(1,0,0,0),(1,1,0,0),(\1,\1,\1,\1),(\1,\1,\1,-\1)\}$

$\ell:=2z_1$

$W:=(\ZZ/2\ZZ)^2=\{(z_1,z_2,z_3,z_4)\mapsto(z_1,\pm z_2,\pm z_3,z_4)\}$

$\Delta^+=\{2z_2,2z_3\}$

$\Phi^+=\{z_1\pm z_2,z_2\pm z_3,z_3\pm z_4\}$

$\pm W$-orbits of $\Sigma$: $[R,\pm S,T,-T]$

$\rho=(r+s+t,s+t,t,0)$}

\noindent
$Sp_{2p}(\CC)\times\CC^*\times\CC^*$ on $\CC^{2p}\oplus\CC^{2p}$ with $2\le
p$

$r=\1$, $s=p-1$, $t=\1$, $e_1=\epsilon+\epsilon'$, $e_2=\omega_2$,
$e_3=\alpha_1$, $e_4=\epsilon-\epsilon'$

\beginrefs

\L|Abk:BenRat1|Sig:\\BR|Au:Benson, C., Ratcliff, G.|Tit:A
classification of multiplicity free actions%
|Zs:J. Algebra|Bd:181|S:152--186|J:1996||

\L|Abk:BenRat2|Sig:\\BR|Au:Benson, C.; Ratcliff, G.
|Tit:Combinatorics and spherical functions on the Heisenberg group%
|Zs:Represent. Theory {\rm(electronic)}|Bd:2|S:79--105|J:1998||

\B|Abk:Bou|Sig:Bou|Au:Bourbaki, N.|Tit:Groupes et alg\`ebres des %
Lie. Chap.~4,~5,~6|Reihe:-|Verlag:Masson|Ort:Paris|J:1981||

\Pr|Abk:He|Sig:He|Au:Heckman, G.|Artikel:Hypergeometric and spherical %
functions|Titel:Harmonic functions and spherical functions on %
symmetric spaces|Hgr:-|Reihe:Perspectives in Mathematics%
|Bd:16|Verlag:Academic Press|Ort:San Diego|S:1--89|J:1994||

\L|Abk:HoUm|Sig:HU|Au:Howe, R., Umeda, T.|Tit:The Capelli
identity, the double commutant theorem, and multiplicity-free
actions|Zs:Math. Ann.|Bd:290|S:565--619|J:1991||

\L|Abk:Kac|Sig:Kac|Au:Kac, V.|Tit:Some remarks on nilpotent orbits%
|Zs:J. Algebra|Bd:64|S:190--213|J:1980||

\Pr|Abk:Montreal|Sig:\\Kn|Au:Knop, F.|Artikel:Some remarks on %
multiplicity free spaces|Titel:Proc. NATO Adv. Study Inst. on %
Representation Theory and Algebraic Geometry|Hgr:A.~Broer, %
G.~Sabi\-dussi, eds.|Reihe:Nato ASI Series C|Bd:514
|Verlag:Kluwer|Ort:Dortrecht|S:301--317|J:1998||

\L|Abk:Semisym|Sig:\\Kn|Au:Knop, F.|Tit:Semisymmetric polynomials %
and the invariant theory of matrix vector pairs|Zs:Preprint %
\tt math.RT/9910060|Bd:-|S:26 pages|J:1999||

\L|Abk:SymCap|Sig:KS|Au:Knop, F.; Sahi, S.|Tit:Difference %
equations and symmetric polynomials defined by their
zeros|Zs:Internat. Math. Res. Notices|Bd:10|S:473--486|J:1996||

\L|Abk:Leahy|Sig:Le|Au:Leahy, A.|Tit:A classification of
multiplicity free representations|Zs:J. Lie Theory|Bd:8|S:367--391|J:1998||

\L|Abk:Ok2|Sig:Ok|Au:Okounkov, A.|Tit:${\rm BC}$-type interpolation Macdonald
polynomials and binomial formula for Koornwinder
polynomials|Zs:Transform. Groups|Bd:3|S:181--207|J:1998||

\L|Abk:OO1|Sig:\\OO|Au:Olshanski, G.; Okounkov, A.|Tit:Shifted Schur
functions|Zs:St.~Petersburg Math. J.|Bd:9|S:73--146|J:1998||

\L|Abk:OO2|Sig:\\OO|Au:Olshanski, G.; Okounkov, A.|Tit:Shifted Jack
polynomials, binomial formula, and applications|Zs:Math. Res. Lett.%
|Bd:4|S:69--78|J:1997||

\L|Abk:ViKi|Sig:VK|Au:Vinberg, E., Kimelfeld,
B.|Tit:Homogeneous domains on flag manifolds and spherical subsets of
semisimple Lie groups|Zs:Funktsional. Anal. i
Prilozhen.|Bd:12|S:12--19|J:1978||

\L|Abk:Yan|Sig:Yan|Au:Yan, Zhi Min|Tit:Special functions associated
with multiplicity free representations|Zs:Preprint|Bd:-|S:-|J:-||

\endrefs

\bye